%% file: MPI_LLG_arxiv.tex
\documentclass[final,leqno]{siamltex}
\usepackage{chngcntr}

\usepackage[arrow, matrix, curve]{xy}
\usepackage[latin1]{inputenc} 	
\usepackage[T1]{fontenc}
\usepackage[english]{babel}
\usepackage{lmodern}
\usepackage{amsmath}
\usepackage{amssymb}  
\usepackage{mathrsfs}
\usepackage{euscript}
\usepackage{enumerate}
\usepackage{xspace}
\usepackage{hyperref}  
\usepackage{framed}
\usepackage{cite}
\usepackage{fancyhdr}
\usepackage{bbm}
\usepackage{float}
\usepackage{comment}
\usepackage{graphicx}
\usepackage{textcomp}
\usepackage{url}
\usepackage{srcltx} 
\usepackage{hyperref} 

\DeclareFontShape{OMX}{cmex}{m}{n}{
  <-7.5> cmex7
  <7.5-8.5> cmex8
  <8.5-9.5> cmex9
  <9.5-> cmex10
}{}
\SetSymbolFont{largesymbols}{normal}{OMX}{cmex}{m}{n}
\SetSymbolFont{largesymbols}{bold}  {OMX}{cmex}{m}{n}

\hypersetup{draft=false, colorlinks=true, linkcolor=black, urlcolor=blue, citecolor=black}

\newcommand{\dd}{\, \mathrm{d}}
\newcommand{\RR}{\mathbb{R}}

\newcommand{\NN}{\mathbb{N}}

\newcommand{\fx}{\mathbf{x}}
\newcommand{\fv}{\mathbf{v}}
\newcommand{\fa}{\mathbf{a}}
\newcommand{\fb}{\mathbf{b}}
\newcommand{\fc}{\mathbf{c}}
\newcommand{\fM}{\mathbf{M}}
\newcommand{\fB}{\mathbf{B}}
\newcommand{\fA}{\mathbf{A}}
\newcommand{\fH}{\mathbf{H}}
\newcommand{\fh}{\mathbf{h}}
\newcommand{\fm}{\mathbf{m}}
\newcommand{\pR}{\mathbf{p}^{\mathrm{R}}}
\newcommand{\Heff}{\fH_{\mathrm{eff}}}
\newcommand{\Hext}{\fH_{\mathrm{ext}}}
\newcommand{\hext}{\fh_{\mathrm{ext}}}

\newcommand{\laplace}{\Delta}
\newcommand{\cX}{\mathcal{X}}
\newcommand{\cY}{\mathcal{Y}}

\def\R{\mathbb{R}}
\def\bF{\mathbb{F}}
\def\altil{\hat{\alpha}}
\def\mv{\mathbf{m}}
\def\mh{\hat{\mathbf{m}}}
\def\uv{\mathbf{u}}
\def\vv{\mathbf{v}}
\def\yv{\mathbf{y}}
\def\wv{\mathbf{w}}
\def\zv{\mathbf{z}}
\def\hv{\mathbf{h}}
\def\Tt{T}
\def\ker{\mathbf{K}}
\def\cU{\mathcal{U}}

  \newtheorem{theo}{Theorem}[section]

  \newtheorem{remark}[theo] {Remark}

\def\cK{\mathcal{K}}
\def\cX{\mathcal{X}}
\def\cY{\mathcal{Y}}
\def\cB{\mathcal{B}}
\def\cV{\mathcal{V}}
\def\cW{\mathcal{W}}
\def\cUtil{\tilde{\mathcal{U}}}
\def\n{\textbf{n}}
\def\m{\textbf{m}}
\def\M{\textbf{M}}
\def\u{\textbf{u}}
\def\vtil{\mathbf{\tilde{v}}}
\def\ve{\mathbf{v^\epsilon}}
\def\vet{\mathbf{v_t^\epsilon}}
\def\be{\mathbf{b^\epsilon}}
\def\pp{\textbf{q}}
\def\p{\pp^z}
\def\pr{\textbf{p}_{\ell}^R}

\def\h{\textbf{h}}

\def\Ktil{\tilde{K}}
\def\KtilT{\tilde{K}_T}

\def\Om2R{\Omega,\R^3}
\def\Pr{\mbox{Proj}_{\m^\bot}}
\def\Prn{\mbox{Proj}_{\n^\bot}}

\author{Barbara Kaltenbacher\thanks{Alpen-Adria-Universit\"at Klagenfurt, Austria ({\tt barbara.kaltenbacher@aau.at})} \and Tram Thi Ngoc Nguyen\thanks{Alpen-Adria-Universit\"at Klagenfurt, Austria ({\tt tram.nguyen@aau.at})} \and Anne Wald\thanks{Department of Mathematics, Saarland University, PO Box 15 11 50, 66123 Saarbr\"ucken, Germany ({\tt anne.wald@num.uni-sb.de})} \and Thomas Schuster\thanks{Department of Mathematics, Saarland University, PO Box 15 11 50, 66123 Saarbr\"ucken, Germany ({\tt thomas.schuster@num.uni-sb.de})}}

\title{Parameter identification for the Landau-Lifshitz-Gilbert equation in Magnetic Particle Imaging}
        
\begin{document}

\maketitle

\begin{abstract}
Magnetic particle imaging (MPI) is a tracer-based technique for medical imaging where the tracer consists of ironoxide nanoparticles. The key idea is to measure the particle response to a temporally changing external magnetic field to compute the spatial concentration of the tracer inside the object. A decent mathematical model demands for a data-driven computation of the system function which does not only describe the measurement geometry but also encodes the interaction of the particles with the external magnetic field. The physical model of this interaction is given by the Landau-Lifshitz-Gilbert (LLG) equation. The determination of the system function can be seen as an inverse problem of its own which can be interpreted as a calibration problem for MPI. In this contribution the calibration problem is formulated as an inverse parameter identification problem for the LLG equation. We give a detailed analysis of the direct as well as the inverse problem in an all-at-once as well as in a reduced setting. The analytical results yield a deeper understanding of inverse problems connected to the LLG equation and provide a starting point for the development of robust numerical solution methods in MPI.
\end{abstract}


\input{introduction}

\input{model}

\input{inverse}

\input{diffadj}

\input{conclusion}

\vspace*{2ex}

\textbf{Acknowledgements:} The work of Anne Wald and Thomas Schuster was partly funded by Hermann und Dr. Charlotte Deutsch-Stiftung and by the German Federal Ministry of Education and Research (Bundesministerium f\"ur Bildung und Forschung, BMBF) under 05M16TSA. 

\bibliographystyle{siam}
\bibliography{bibliography_llg}

\end{document}

%% file: introduction.tex
\section{Introduction}
Magnetic particle imaging (MPI) is a dynamic imaging modality for medical applications that has first been introduced in 2005 by B.~Gleich and J.~Weizenecker \cite{bgjw05}. Magnetic nanoparticles, consisting of a magnetic iron oxide core and a nonmagnetic coating, are inserted into the body to serve as a tracer. The key idea is to measure the nonlinear response of the nanoparticles to a temporally changing external magnetic field in order to draw conclusions on the spatial concentration of the particles inside the body. Since the particles are distributed along the bloodstream of a patient, the particle concentration yields information on the blood flow and is thus suitable for cardiovascular diagnosis or cancer detection \cite{tktb12, tkngmm17}. An overview of MPI basics is given in \cite{tktb12}. Since MPI requires the nanoparticles as a tracer, it mostly yields quantitative information on their distribution, but does not image the morphology of the body, such as the tissue density. The latter can be visualized using computerized tomography (CT) \cite{natterer86} or magnetic resonance imaging (MRI) \cite{HINSHAW;LENT:83}. These do not require a tracer, but involve ionizing radiation in the case of CT or, in the case of MRI, a strong magnetic field and a potentially high acquisition time. Other tracer-based methods are, e.g., single photon emission computerized tomography (SPECT) and positron emission tomography (PET) \cite{natterer_buch_2, SHEPP;VARDI:82}, which both involve radioactive radiation. The magnetic nanoparticles that are used in MPI, on the other hand, are not harmful for organisms. For a more detailed comparison of these methods, we would like to refer the reader to \cite{tktb12}.   

At this point there have been promising preclinical studies on the performance of MPI, showing that this imaging modality has a great potential for medical diagnosis since it is highly sensitive with a good spatial and temporal resolution, and the data acquisition is very fast \cite{tkngmm17}. However, particularly in view of an application to image the human body, there remain some obstacles. One obstacle is the time-consuming calibration process. In this work, we assume that the concentration of the nanoparticles inside the body remains static throughout both the calibration process and the actual image acquisition. Mathematically, the forward problem of MPI then can essentially be formulated as an integral equation of the first kind for the particle concentration (or distribution) $c$, 
\begin{displaymath}
 u(t) = \int_{\Omega} c(x) s(x,t) \dd x,
\end{displaymath}
where the integration kernel $s$ is called the \emph{system function}. The system function encodes some geometrical aspects of the MPI scanner, such as the coil sensitivities of the receive coils in which the particle signal $u$ is measured, but mostly it is determined by the particle behavior in response to the applied external magnetic field.  

The actual inverse problem in MPI is to reconstruct the concentration $c$ under the knowledge of the system function $s$ from the measured data $u$. To this end, the system function has to be determined prior to the scanning procedure. This is usually done by evaluating a series of full scans of the field of view, where in each scan a delta sample is placed in a different pixel until the entire field of view is covered \cite{tktb12}. Another option is a model-based approach for $s$ (see for example \cite{tkpm17}), which basically involves a model for the particle magnetization. Since this model often depends on unknown parameters, the model-based determination of the system function itself can again be formulated as an inverse problem. This article now addresses this latter type of inverse problem, i.e., the identification of the system function for a known set of concentrations from calibration measurements. More precisely, our goal is to find a decent model for the time-derivative of the particle magnetization $\fm$, which is proportional to $s$.

So far, in model-based approaches for the system function, the particle magnetization $\fm$ is not modeled directly. Instead, one describes the mean magnetization $\overline{\fm}$ of the particles via the Langevin function, i.e., the response of the particles is modeled on the mesoscopic scale \cite{tktb12, tk18}. This approach is based on the assumption that the particles are in thermodynamic equilibrium and respond directly to the external field. For this reason, the mean magnetization is assumed to be a function of the external field, such that the mean magnetization is always aligned with the external field. The momentum of the mean magnetization is calculated via the Langevin function. This model, however, neglects some properties of the particle behavior. In particular, the magnetic moments of the particles do not align instantly with the external field \cite{cgc12}. 

In this work, we thus address an approach from micromagnetics, which models the time-dependent behavior of the magnetic material inside the particles' cores on the micro scale and allows to take into account various additional physical properties such as particle-particle interaction. For an overview, see for example \cite{mkap06}. Since the core material is iron oxide, which is a ferrimagnetic material that shows a similar behavior as ferromagnets \cite{cg11imm, demt2}, we use the \emph{Landau-Lifshitz-Gilbert (LLG) equation}
\begin{displaymath}
 \frac{\partial}{\partial t} \fm = - \widetilde{\alpha}_1 \fm \times \left( \fm \times \fH_{\mathrm{eff}} \right) + \widetilde{\alpha}_2 \fm \times \Heff,
\end{displaymath}
see also \cite{tlg04, llel92}, for the evolution of the magnetization $\fm$ of the core material. The field $\fH_{\mathrm{eff}}$ incorporates the external magnetic field together with other relevant physical effects. According to the LLG equation, the magnetization $\fm$ performs a damped precession around the field vector of the external field, which leads to a relaxation effect. The LLG equation has been widely applied to describe the time evolution in micromagnetics \cite{bsffgppr14, bgmh93, dkps15}.

In contrast to the imaging problem of MPI, the inverse problem of determining the magnetization $\fm$ along with the constants $\widetilde{\alpha}_1, \widetilde{\alpha}_2$ turns out to be a nonlinear inverse problem, which is typical for parameter identification problems for partial differential equations, for example electrical impedance tomography \cite{borcea02}, terahertz tomography \cite{awts17}, ultrasound imaging \cite{ck13} and other applications from imaging and nondestructive testing \cite{akar16}.\\
We use the \emph{all-at-once} as well as the \emph{reduced} formulation of this inverse problem in a Hilbert space setting, see also \cite{bk16, bk17, ttnn19}, and analyze both cases including well-definedness of the forward mapping, continuity, and Fr\'echet differentiability and calculate the adjoint mappings for the Fr\'echet derivatives. By consequence, iterative methods such as the Landweber method \cite{HaNeSc95, Landweber}, also in combination with Kaczmarz' method \cite{HKLS07,HLS07}, Newton methods (see, e.g., \cite{rieder99}), or subspace techniques \cite{ws16} can be applied for the numerical solution. An overview of suitable regularization techniques is given in \cite{kns08, kss09}.

We begin with a detailed introduction to the modelling in MPI. In particular, we describe the full forward problem and present the initial boundary value problem for the LLG equation that we use to describe the magnetization evolution. In Section \ref{sec_inverse}, we formulate the inverse problem of calibration both in the all-at-once and in the reduced setting to obtain the final operator equation that is analyzed in the subsequent section. First, in Section \ref{sec:aao}, we present an analysis for the all-at-once setting. The inverse problem in the reduced setting is then addressed in Section \ref{sec:red}. Finally, we conclude our findings in Section \ref{sec_concl} and give an outlook on further research.

%

\vspace*{2ex}

Throughout the article, we make use of the following notation: The differential operators $-\Delta$ and $\nabla$ are applied by components to a vector field. 
In particular this means that by $\nabla\uv$ we denote the transpose of the Jacobian of $\uv$.
Moreover, $\langle\mathbf{a},\mathbf{b}\rangle$ or $\mathbf{a}\cdot\mathbf{b}$ denotes the Euclidean inner product between two vectors and $A\colon B$ the Frobenius inner product between two matrices.

%% file: model.tex
\section{The underlying physical model for MPI}
The basic physical principle that is exploited in MPI is Faraday's law of induction, which states that whenever the magnetic flux density $\mathbf{B}$ through a coil changes in time, this change induces an electric current in the coil. This current, or rather the respective voltage, can be measured. In MPI, the magnetic flux density $\fB$ consists of the external applied magnetic field $\Hext$ and the particle magnetization $\fM^{\mathrm{P}}$, i.e.,
\begin{displaymath}
 \fB = \mu_0 \left( \Hext + \fM^{\mathrm{P}} \right),
\end{displaymath}
where $\mu_0$ is the magnetic permeability in vacuum. The particle magnetization $\fM^{\mathrm{P}}(x,t)$ in $x \in \Omega \subseteq \RR^3$ depends linearly on the concentration $c(x)$ of magnetic material, which corresponds to the particle concentration, in $x \in \Omega$ and on the magnetization $\fm(x,t)$ of the magnetic material. We thus have
\begin{displaymath}
 \fM^{\mathrm{P}}(x,t) = c(x) \fm(x,t), 
\end{displaymath}
where $\lvert \fm \rvert = m_{\mathrm{S}} > 0$, i.e., the vector $\fm$ has the fixed length $m_{\mathrm{S}}$ that depends on the magnetic core material inside the particles. At this point it is important to remark that we use a slightly different approach to separate the particle concentration, which carries the spatial information on the particles, from the magnetization behavior of the magnetic material and the measuring process. In our approach, the concentration is a dimensionless quantity, whereas in most models, it is defined as the number of particles per unit volume (see, e.g. \cite{tktb12}). 

\vspace*{2ex}

A detailed derivation of the forward model in MPI, based on the equilibrium model for the magnetization, can be found in \cite{tktb12}. The steps that are related to the measuring process can be adapted to our approach. For the reader's convenience, we want to give a short overview and introduce the parameters related to the scanner setup. \\
If the receive coil is a simple conductor loop, which encloses a surface $S$, the voltage that is induced can be expressed by
\begin{equation} \label{faraday}
 u(t) = -\frac{\mathrm{d}}{\mathrm{d}t} \int_{S} \fB(x,t) \cdot \dd \fA = -\mu_0 \frac{\mathrm{d}}{\mathrm{d}t} \int_{S} \left( \Hext + \fM^{\mathrm{P}} \right) \cdot \dd \fA.
\end{equation}
The signal that is recorded in the receive coil thus originates from temporal changes of the external magnetic field $\fH$ as well as of the particle magnetization $\fM^{\mathrm{P}}$,
\begin{align}
 u(t) &=  -\mu_0 \left( \int_{\Omega} \pR(x) \cdot \frac{\partial}{\partial t} \Hext(x,t) \dd x + \int_{\Omega} \pR(x) \cdot \frac{\partial}{\partial t} \fM^{\mathrm{P}}(x,t) \dd x \right) \\
      &=: u^{\mathrm{E}}(t) + u^{\mathrm{P}}(t)
\end{align}
For the signal that is caused by the change in the particle magnetization we obtain
\begin{displaymath}
 \begin{split}
 u^{\mathrm{P}}(t) &= -\mu_0 \frac{\mathrm{d}}{\mathrm{d}t} \int_{\Omega} \pR(x) \cdot \fM^{\mathrm{P}}(x,t) \dd x \\
                   &= -\mu_0 \int_{\Omega} \pR(x) \cdot \frac{\partial}{\partial t} \fM^{\mathrm{P}}(x,t) \dd x \\
                   &= -\mu_0 \int_{\Omega} c(x) \pR(x) \cdot \frac{\partial}{\partial t} \fm(x,t) \dd x \\
                   &= -\mu_0 \int_{\Omega} c(x) s(x,t) \dd x.
 \end{split}
\end{displaymath}
The function 
\begin{equation}\label{systemfunction}
 s(x,t) := \pR(x) \cdot \frac{\partial}{\partial t} \fm(x,t) = \left\langle \pR(x), \frac{\partial}{\partial t} \fm(x,t) \right\rangle_{\RR^3}
\end{equation}
is called the \emph{system function} and can be interpreted as a potential to induce a signal in the receive coil. The function $\pR$ is called the coil sensitivity and is determined by the architecture of the respective receive coil. For our purposes, we assume that $\pR$ is known. The measured signal that originates from the magnetic particles can thus essentially be calculated via an integral equation of the first kind with a time-dependent integration kernel $s$.

\vspace*{2ex}

The particle magnetization, however, changes in time in response to changes of the external field. It is thus an important objective to encode the interplay of the external field and the particles in a sufficiently accurate physical model. The magnetization of the magnetic particles that are used in MPI can be considered on different scales. The following characterization from ferromagnetism has been taken from \cite{mkap06}: \\
On the \emph{atomic level}, one can describe the behavior of a magnetic material as a spin system and take into account stochastic effects that arise, for example, from Brownian motion. \\
In the \emph{microscopic scale}, continuum physics is applied to work with deterministic equations describing the magnetization of the magnetic material. \\
In the \emph{mesoscopic scale}, we can describe the magnetization behavior via a mean magnetization, which is an average particle magnetic moment. \\
Finally, on a \emph{macroscopic scale}, all aspects that arise from the microstructure are neglected and the magnetization is described by phenomenological constitutive laws.

\vspace*{2ex}

In this work, we intend to use a model from micromagnetism, allowing us to work with a deterministic equation to describe the magnetization of the magnetic material. The core material of the nanoparticles consists of iron-oxide or magnetite, which is a ferrimagnetic material. The magnetization curve of ferrimagnetic materials is similar to the curve that is observed for ferromagnets, but with a lower saturation magnetization (see, e.g., \cite{cg11imm, demt2}). This approach has also been suggested in \cite{drjb14}. 
The evolution of the magnetization in time is described by the \emph{Landau-Lifshitz-Gilbert (LLG) equation}
\begin{equation}
 \fm_t := \frac{\partial}{\partial t} \fm = - \widetilde{\alpha}_1 \fm \times \left( \fm \times \fH_{\mathrm{eff}} \right) + \widetilde{\alpha}_2 \fm \times \Heff,
\end{equation}
see \cite{tlg04, mkap06} and the therein cited literature. The coefficients
\begin{displaymath}
 \widetilde{\alpha}_1 := \frac{\gamma \alpha_{\mathrm{D}}}{m_{\mathrm{S}}(1+\alpha_{\mathrm{D}}^2)} > 0, \quad \widetilde{\alpha}_2 := \frac{\gamma}{(1+\alpha_{\mathrm{D}}^2)} > 0
\end{displaymath}
are material parameters that contain the gyromagnetic constant $\gamma$, the saturation magnetization $m_{\mathrm{S}}$ of the core material and a damping parameter $\alpha_{\mathrm{D}}$. The vector field $\Heff$ is called the \emph{effective magnetic field}. It is defined as the negative gradient $-D\mathcal{E}(\fm)$ of the \emph{Landau energy} $\mathcal{E}(\fm)$ of a ferromagnet, see, e.g., \cite{mkap06}. Taking into account only the interaction with the external magnetic field $\fH$ and particle-particle interactions, this energy is given by
\begin{displaymath}
 \mathcal{E}_{A}(\fm) = A \int_{\Omega} \lvert \nabla \fm \rvert^2 \dd x - \mu_0 m_{\mathrm{S}} \int_{\Omega} \left\langle \fH, \fm \right\rangle_{\RR^3} \dd x,
\end{displaymath}
where $A \geq 0$ is a scalar parameter (the exchange stiffness constant \cite{tlg04}). We thus have
\begin{equation}
 \Heff = 2A \Delta \fm + \mu_0 m_{\mathrm{S}}\Hext.
\end{equation}
Together with Neumann boundary conditions and a suitable initial condition our model for the magnetization thus reads
\begin{align}
 \fm_t &= - \alpha_1 \fm \times \left( \fm \times (\Delta \fm + \hext) \right) + \alpha_2 \fm \times (\Delta \fm + \hext) \ && \text{in} \ [0,T] \times \Omega, \label{llg1}\\
 0 &= \partial_{\nu} \fm  && \text{on} \ [0,T] \times \partial\Omega, \label{llg1_bc}\\
 \fm_0 &= \fm(t=0), \ \lvert \fm_0 \rvert = m_{\mathrm{S}} && \text{in} \ \Omega, \label{llg1_abs}
\end{align}
where $\hext = \frac{\mu_0 m_{\mathrm{S}}}{2A} \Hext$ and $\alpha_1 := 2A\widetilde{\alpha}_1, \alpha_2 := 2A\widetilde{\alpha}_2 > 0$. The initial value $\fm_0 = \fm(t=0)$ corresponds to the magnetization of the magnetic material in the beginning of the measurement. To obtain a reasonable value for $\fm_0$, we take into account that the external magnetic field is switched on before the measuring process starts, i.e., $\fm_0$ is the state of the magnetization that is acquired when the external field is static. This allows us to precompute $\fm_0$ as the solution of the stationary problem
\begin{equation} \label{llg_initial_value}
 \alpha_1 \fm_0 \times \left( \fm_0 \times (\Delta \fm_0 + \hext(t=0)) \right) = \alpha_2 \fm_0 \times (\Delta \fm_0 + \hext(t=0))
\end{equation}
with Neumann boundary conditions. 

\begin{remark}
In the stationary case, damping does not play a role, and if we additionally neglect particle-particle interactions, we obtain the approximative equation
\begin{displaymath}
 \hat{\fm}_0 \times \left( \hat{\fm}_0 \times \hext(t=0) \right) = 0
\end{displaymath}
with an approximation $\hat{\fm}_0$ to $\hat{\fm}$, since $\alpha_2 \approx 0$ and $\Heff \approx \mu_0 m_{\mathrm{S}}\Hext$. The above equation yields $\hat{\fm}_0 \parallel \hext(t=0)$. Together with $\lvert \hat{\fm}_0 \rvert = m_{\mathrm{S}}$ this yields
\begin{displaymath}
 \hat{\fm}_0 = m_{\mathrm{S}} \frac{\hext(t=0)}{\lvert \hext(t=0) \rvert}.
\end{displaymath}
This represents a good approximation to $\fm_0$ where $\hext$ is strong at the time point $t=0$:
\begin{displaymath}
 \fm_0 \approx \hat{\fm}_0 = m_{\mathrm{S}} \frac{\hext(t=0)}{\lvert \hext(t=0) \rvert}.
\end{displaymath}
\end{remark}


\subsection{The observation operator in MPI}
Faraday's law states that a temporally changing magnetic field induces an electric current in a conductor loop or coil, which yields the relation \eqref{faraday}. By consequence, not only the change in the particle magnetization contributes to the induced current, but also the dynamic external magnetic field $\Hext$. Since we need the particle signal for the determination of the particle magnetization, we need to separate the particle signal from the excitation signal due to the external field. This is realized by processing the signal in a suitable way using filters. \\
MPI scanners usually use multiple receive coils to measure the induced particle signal at different positions in the scanner. We assume that we have $L \in \NN$ receive coils with coil sensitivities $\pR_\ell$, $\ell=1,...,L$, and the measured signal is given by
\begin{equation}
 \widetilde{v}_\ell(t) = -\mu_0 \int_{0}^T \widetilde{a}_\ell(t - \tau ) \int_{\Omega} c(x) \pR_\ell(x) \cdot \frac{\partial}{\partial \tau} \fm(x,\tau) \dd x \dd \tau,
\end{equation}
where $T$ is the repetition time of the acquisition process, i.e., the time that is needed for one full scan of the object, and $a_\ell \, : \, [0,T] \rightarrow \RR$ is the transfer function with periodic continuation $\widetilde{a}_\ell \, : \, \RR \rightarrow \RR$. The transfer function serves as a filter to separate particle and excitation signal, i.e., it is chosen such that
\begin{displaymath}
 \widetilde{v}^{\mathrm{E}}_\ell(t) := \big(\widetilde{a}_\ell * u^{\mathrm{E}}_\ell\big)(t) = -\mu_0\int_0^T \widetilde{a}_\ell(t-\tau)\int_{\Omega} \pR_\ell(x) \cdot \frac{\partial}{\partial t} \Hext(x,t) \dd x \dd t \approx 0.
\end{displaymath}
In practice, $\widetilde{a}_\ell$ is often a band pass filter. For a more detailed discussion of the transfer function, see also \cite{tktb12}. In this work, the transfer function is known analytically.

We define
\begin{displaymath}
 \mathbf{K}_\ell(t,\tau,x) := -\mu_0\widetilde{a}_\ell(t - \tau ) c(x) \pR_\ell(x),
\end{displaymath}
such that the measured particle signals are given by
\begin{equation} \label{forwardoperator}
 {v}_\ell(t) = \int_{0}^T \int_{\Omega} \mathbf{K}_\ell(t,\tau,x) \cdot \frac{\partial}{\partial \tau} \fm(x,\tau) \dd \tau \dd x,
\end{equation}
where $\fm$ fulfills \eqref{llg1}, \eqref{llg1_bc}, \eqref{llg1_abs}.

\vspace*{2ex}

To determine $\fm$ in $\Omega\times(0,T)$, we use the data ${v}_{k\ell}(t)$, $k=1,...,K$, $\ell=1,...,L$, from the scans that we obtain for different particle concentrations $c_k$, $k=1,...,K$, $K \in \NN$. The forward operator thus reads
\begin{equation} \label{forwardoperator_kl}
\begin{split}
& {v}_{k\ell}(t) = \int_{0}^T \int_{\Omega} \mathbf{K}_{k\ell}(t,\tau,x) \cdot \frac{\partial}{\partial \tau}\fm(x,\tau) \dd x \dd \tau \,, \\
&\mathbf{K}_{k\ell}(t,\tau,x) := -\mu_0\widetilde{a}_\ell(t - \tau ) c_k(x) \pR_\ell(x).
\end{split}
\end{equation}

\subsection{Equivalent formulations of the LLG equation}
In this section, we derive additional formulations of \eqref{llg1} - \eqref{llg1_abs} that are suitable for the analysis. The approach is motivated by \cite{mkap06}, where only particle-particle interactions are taken into account. \\ 
First of all, we observe that multiplying \eqref{llg1} with $\fm$ on both sides yields
\begin{equation}
 \frac{1}{2} \cdot \frac{\mathrm{d}}{\mathrm{d}t} \lvert \fm(x,t) \rvert^2 = \fm(x,t) \cdot \fm_t(x,t) = 0,
\end{equation}
which shows that the absolute value of $\fm$ does not change in time. Since $\lvert \fm_0 \rvert = m_{\mathrm{S}}$, we have $\fm(x,t) \in m_{\mathrm{S}}\cdot \mathcal{S}^2$, where $\mathcal{S}^2 := \lbrace \fv \in \RR^3 \, : \, \lvert \fv \rvert = 1 \rbrace$ is the unit sphere in $\RR^3$. 
As a consequence, we have $0 = \nabla \lvert \fm \rvert^2 = 2\nabla \fm \cdot \fm$ in $\Omega$, so that, by taking the divergence we get
\begin{equation} \label{grad_abs_m}
 \langle \fm, \laplace \fm \rangle  = - \langle \nabla \fm, \nabla \fm \rangle.
\end{equation}
Now we make use of the identity
\begin{displaymath}
 \fa \times (\fb \times \fc) = \langle \fa,\fc \rangle \fb - \langle \fa,\fb \rangle \fc
\end{displaymath}
for $\fa,\fb,\fc \in \RR^3$ to derive
\begin{align}
 \fm \times (\fm \times \laplace \fm) &= \langle \fm, \laplace \fm \rangle \fm - \lvert \fm \rvert^2 \laplace \fm
= -\lvert \nabla \fm \rvert^2\fm-m_{\mathrm{S}}^2\laplace \fm, 
\label{bac_cab_1} \\
 \fm \times (\fm \times \hext) &= \langle \fm, \hext \rangle \fm - \lvert \fm \rvert^2 \hext
=\langle \fm, \hext \rangle \fm - m_{\mathrm{S}}^2\hext
. \label{bac_cab_2}
\end{align}
Using \eqref{grad_abs_m} together with \eqref{bac_cab_1}, \eqref{bac_cab_2} and $\lvert \fm \rvert = m_{\mathrm{S}}$, we obtain from \eqref{llg1} - \eqref{llg1_abs}
\begin{align}
\begin{split}
 \fm_t - \alpha_1 \, m_{\mathrm{S}}^2\,\laplace \fm &= \alpha_1 \lvert \nabla \fm \rvert^2 \fm + \alpha_2 \fm \times \laplace \fm  \\
 &\quad- \alpha_1 \langle \fm, \hext \rangle \fm + \alpha_1 \,m_{\mathrm{S}}^2\,\hext + \alpha_2 \fm \times \hext 
\end{split} && \text{in} \ [0,T] \times \Omega, \label{llg2}\\
 0 &= \partial_{\nu} \fm  && \text{on} \ [0,T] \times \partial\Omega, \label{llg2_bc}\\
 \fm_0 &= \fm(t=0), \ \lvert \fm_0 \rvert = m_{\mathrm{S}} && \text{in} \ \Omega, \label{llg2_abs}
\end{align}

Taking the cross product of $\fm$ with \eqref{llg2} and multiplying with $-\hat{\alpha}_2$, where $\hat{\alpha}_1=\frac{\alpha_1}{m_{\mathrm{S}}^2\alpha_1^2+\alpha_2^2}$, $\hat{\alpha}_2=\frac{\alpha_2}{m_{\mathrm{S}}^2\alpha_1^2+\alpha_2^2}$, by \eqref{bac_cab_1}, \eqref{bac_cab_2} and cancellation of the first and third term on the right hand side we get
\[
\begin{aligned}
& -\hat{\alpha}_2\fm\times\fm_t + \alpha_1\hat{\alpha}_2m_{\mathrm{S}}^2\,\fm\times\laplace \fm \\
&=  
\frac{\alpha_2^2}{m_{\mathrm{S}}^2\alpha_1^2+\alpha_2^2} \left(\lvert \nabla \fm \rvert^2\fm+m_{\mathrm{S}}^2\laplace \fm\right)\\  
&\hspace{2cm} - \alpha_1\hat{\alpha}_2 m_{\mathrm{S}}^2\,\fm\times\hext 
+ \frac{\alpha_2^2}{m_{\mathrm{S}}^2\alpha_1^2+\alpha_2^2} \left(m_{\mathrm{S}}^2\hext - \langle \fm, \hext \rangle \fm\right)\,,
\end{aligned}
\]
where the second term on the left hand side can be expressed via \eqref{llg2} as 
\[
\alpha_1\hat{\alpha}_2\fm\times\laplace \fm =  \hat{\alpha}_1\fm_t + \frac{\alpha_1^2}{m_{\mathrm{S}}^2\alpha_1^2+\alpha_2^2}\left( -m_{\mathrm{S}}^2\laplace \fm - \lvert \nabla \fm \rvert^2 \fm  + \langle \fm, \hext \rangle \fm - m_{\mathrm{S}}^2\hext\right) 
-\alpha_1\hat{\alpha}_2 \fm \times \hext\,.
\]
This yields the alternative formulation
\begin{align}
\hat{\alpha}_1 m_{\mathrm{S}}^2 \fm_t - \hat{\alpha}_2 \fm \times \fm_t -m_{\mathrm{S}}^2\laplace \fm &= \lvert \nabla \fm \rvert^2 \fm + m_{\mathrm{S}}^2\hext - 
\langle \fm, \hext \rangle \fm  && \text{in} \ [0,T] \times \Omega, \label{llg3}\\
 0 &= \partial_{\nu} \fm  && \text{on} \ [0,T] \times \partial\Omega, \label{llg3_bc}\\
 \fm_0 &= \fm(t=0), \ \lvert \fm_0 \rvert = m_{\mathrm{S}} && \text{in} \ \Omega \label{llg3_abs}\,.
\end{align}

%% file: inverse.tex
\section{An inverse problem for the calibration process in MPI}  \label{sec_inverse}

Apart from the obvious inverse problem of determining the concentration $c$ of magnetic particles inside a body from the measurements $v_{\ell}$, $\ell = 1,...,L$, MPI gives rise to a range of further parameter identification problems of entirely different nature. In this work, we are not addressing the imaging process itself, but consider an inverse problem that is essential for the calibration process. Here, calibration refers to determining the system function $s_{\ell}$, which serves as an integral kernel in the imaging process. The system function includes all system parameters of the tomograph and encodes the physical behaviour of the magnetic material in the cores of the magnetic particles inside a temporally changing external magnetic field. Experiments show that a simple model for the magnetization, based on the assumption that the particles are in their equilibrium state at all times, is insufficient for the imaging, see, e.g., \cite{tkpm17}. A model-based approach with an enhanced physical model has so far been omitted due to the complexity of the involved physics and the system function is usually measured in a time-consuming calibration process \cite{tktb12, tkngmm17}.

In this work, we address the inverse problem of calibrating an MPI system for a given set of standard calibration concentrations $c_k$, $k=1,...,K$, 
for which we measure the corresponding signals and obtain the data ${v}_{k\ell}(t)$, $k=1,...,K$, $\ell=1,...,L$. Here we assume that the coil sensitivity $\pR_{\ell}$ as well as the transfer function $\widetilde{a}_{\ell}$ are known. 

This, together with the fact that $\fm$ is supposed to satisfy the LLG equation \eqref{llg3}--\eqref{llg3_abs}, is used to determine the system function \eqref{systemfunction}. Actually, since $\pR$ is known, the inverse problem under consideration here consists of reconstructing $\fm$ from \eqref{forwardoperator_kl}, \eqref{llg3}--\eqref{llg3_abs}. As the initial boundary value problem \eqref{llg3}--\eqref{llg3_abs} has a unique solution $\fm$ for given $\hat{\alpha}_1$, $\hat{\alpha}_2$, it actually suffices to determine these two parameters. This is the point of view that we take when using a classical reduced formulation of the calibration problem
\begin{equation}\label{ip_red}
F(\altil)=y
\end{equation}
with the data $y_{k\ell}= v_{k\ell}$ and the forward operator
\begin{align} \label{red-F}
F:\mathcal{D}(F)(\subseteq\cX)\rightarrow\cY, \qquad \altil=(\altil_1,\altil_2) \mapsto \mathcal{K} \frac{\partial}{\partial t} S(\altil)
\end{align}
containing the parameter-to-state map 
\begin{align} \label{defS}
S:\cX\rightarrow\cUtil
\end{align}
that maps the parameters $\altil$ into the solution $\m:=S(\altil)$ of the LLG initial boundary value problem \eqref{llg3}--\eqref{llg3_abs}.
The linear operator $\mathcal{K}$ is the integral operator defined by the kernels $\mathbf{K}_{k\ell}$, $k=1,...,K$, $\ell=1,...,L$, i.e., 
\begin{align} \label{defK}
\mathcal{K}_{k\ell} \uv = \int_0^T \int_\Omega \ker_{k\ell}(t,\tau,\fx)\cdot \uv(\fx,\tau)\dd \tau \dd \fx\,.
\end{align}
Here, the preimage and image spaces are defined by  
\begin{align}\label{defXY}
\cX = \R^2, \qquad \cY=L^2(0,T)^{KL}
\end{align}
and the state space $\cUtil$ will be chosen appropriately below, see Section \ref{sec:red}.

Alternatively, we also consider the all-at-once formulation of the inverse problem as a simultaneous system 
\begin{equation}\label{ip_aao}
\bF(\fm,\altil)=\mathbf{y}:=(0,y)^T
\end{equation}
for the state $\fm$ and the parameters $\altil$, with the forward operator
\[
\bF(\fm,\altil)= \left(\begin{array}{l}
\bF_0(\fm,\altil)\\
\Bigl(\bF_{k\ell}(\fm,\altil)\Bigr)_{k=1,\ldots,K\,, \ \ell=1,\ldots,L}
\end{array}\right),
\]
where 
\[
\bF_0(\fm,\altil_1,\altil_2)=:
\altil_1 \fm_t-\Delta \fm - \altil_2 \fm\times \fm_t -|\nabla \fm|^2\fm - \hext +(\fm\cdot\hext)\fm
\]
and 
\[
\bF_{k\ell}(\fm,\altil_1,\altil_2) = \mathcal{K}_{k,\ell} \fm_t
\]
with $\mathcal{K}_{k,\ell}$ as in \eqref{defK}. 
Here $\bF$ maps between $\cU\times\cX$ and $\cW\times\cY$ with $\cX$, $\cY$ as in \eqref{defXY}, and $\cU$, $\cW$ appropriately chosen function spaces, see Section \ref{sec:aao}.

\medskip

Iterative methods for solving inverse problems usually require the linearization $F'(\altil)$ of the forward operator $F$ and its adjoint $F'(\altil)^*$ (and likewise for $\bF$) in the given Hilbert space setting.

For example, consider Landweber's iteration cf., e.g., \cite{Landweber,HaNeSc95} defined by a gradient decent method for the least squares functional $\|F(\altil)-y\|_\cY^2$ as
\[
\altil_{n+1}=\altil_n-\mu_n F'(\altil_n)^*(F(\altil_n)-y)
\]
with an appropriately chosen step size $\mu_n$.
Alternatively, one can split the forward operator into a system by considering it row wise $F_k(\altil)=y_k$ with $F_k=(F_{kl})_{\ell=1\dots L}$ or column wise $F_\ell(\altil)=y_\ell$ with $F_\ell=(F_{kl})_{k=1,\ldots,K}$, or even element wise $F_{kl}(\altil)=y_{kl}$, and cyclically iterating over these equations with gradient descent steps in a Kaczmarz version of the Landweber iteration cf., e.g., \cite{HLS07,HKLS07}.
The same can be done with the respective all-at-once versions \cite{bk16}. These methods extend to Banach spaces as well by using duality mappings, cf., e.g., \cite{SKHK12}, however, for the sake of simplicity of exposition and implementation, we will concentrate on a Hilbert space setting here; in particular, all adjoints will be Hilbert space adjoints.

%% file: diffadj.tex
\section{Derivatives and adjoints}

Motivated by their need in iterative reconstruction methods, we now derive and rigorously justify derivatives of the forward operators as well as their adjoints, both in an all-at-once and in a reduced setting.

To simplify notation for the following analysis sections, the subscript "ext" in the external magnetic field will be skipped. Moreover, to avoid confusion with the dual pairing, we will use the dot notation for the Euclidean inner product.

\input{diffadj_aao}

\input{diffadj_red}

%% file: diffadj_aao.tex
\subsection{All-at-once formulation} \label{sec:aao}
We split the magnetization additively into its given initial value $\mv_0$ and the unknown rest $\mh$, so that the forward operator reads 
\[
\begin{aligned}
&\bF(\mh,\altil_1,\altil_2)=
\left(\begin{array}{l}
\bF_0(\mh,\altil_1,\altil_2)\\[2ex]
\Bigl(\bF_{k\ell}(\mh,\altil_1,\altil_2)\Bigr)_{k=1,\ldots,K\,, \ \ell=1,\ldots,L}
\end{array}\right)\\
&:=
\left(\begin{array}{l}
\altil_1 \mh_t-\Delta_N (\mv_0+\mh) - \altil_2 (\mv_0+\mh)\times \mh_t 
\\ \hspace*{1.5cm}
-|\nabla (\mv_0+\mh)|^2(\mv_0+\mh) - \hv +((\mv_0+\mh)\cdot\hv)(\mv_0+\mh)\\[2ex]
\Bigl(\int_0^\Tt \int_\Omega \ker_{k\ell}(t,\tau,x)\cdot \mv_t(x,\tau)\, dx\, d\tau\Bigr)_{k=1,\ldots,K\,, \ \ell=1,\ldots,L}
\end{array}\right)\,,
\end{aligned}
\]
for given $\hv\in L^2(0,T;L^p(\Omega;\R^3))$, $p\geq2$, where $\Delta_N:H^1(\Omega)\to H^1(\Omega)^*$ and, using the same notation, $\Delta_N:H^2_N(\Omega)\to L^2(\Omega)(\subseteq H^1(\Omega)^*)$ with $H^2_N(\Omega)=\{u\in H^2(\Omega)\, : \, \partial_\nu u=0\mbox{ on }\partial\Omega\}$
\footnote{Note that as opposed to $H^1(\Omega)$ functions, $H^2(\Omega)$ functions do have a Neumann boundary trace}   
is equipped with homogeneous Neumann boundary conditions, i.e, it is defined by
\[
\langle -\Delta_N u,v\rangle_{H^1(\Omega)^*,H^1(\Omega)}= (\nabla u, \nabla v)_{L^2(\Omega)} \quad \forall u,v\in H^1(\Omega) 
\]
and thus satisfies 
\begin{equation}\label{Laplace}
(-\Delta_N u,v)_{L^2(\Omega)}=\int_\Omega \nabla u\cdot \nabla v\, dx \quad \forall u\in H^2_N(\Omega)\,, \ v\in H^1(\Omega)\,.
\end{equation}

The forward operator is supposed to act between Hilbert spaces 
\[
\bF:\cU\times \R^2 \to \cW\times L^2(0,T)^{K L}
\]
with the linear space
\begin{equation}\label{U}
\begin{aligned}
\cU&=\{ \uv\in L^2(0,T;H^2_N(\Omega;\R^3))\cap H^1(0,T;L^2(\Omega;\R^3))\, : \, \uv(0)=0\}\\
&\subseteq C(0,T;H^1(\Omega))\cap H^s(0,T;H^{2-2s}(\Omega))\,,
\end{aligned}
\end{equation}
for $s\in[0,1]$, where the latter embedding is continuous by, e.g, \cite[Lemma 7.3]{Roubicek}, applied to $\frac{\partial u_i}{\partial x_j}$, and interpolation, as well as
\begin{equation}\label{W}
\cW=H^1(0,T;H^1(\Omega;\R^3))^* \mbox{ or, in case $p>2$, } \cW=H^1(0,T;L^2(\Omega;\R^3))^*\,.
\end{equation}
We equip $\cU$ with the inner product
\[
(\uv_1,\uv_2)_\cU:=\int_0^T\int_\Omega \Bigl( (-\Delta_N\uv_1)\cdot(-\Delta_N\uv_2) +\uv_{1t}\cdot\uv_{2t}\Bigr)\, dx\, dt 
+\int_\Omega \nabla\uv_1(T) \colon \nabla\uv_2(T)\, dx\,,
\]
which, in spite of the nontrivial nullspace of the Neumann Laplacian $-\Delta_N$, defines a norm equivalent to the usual norm on $L^2(0,T;H^2(\Omega;\R^3))\cap H^1(0,T;L^2(\Omega;\R^3))$, due to the estimates 
\[
\begin{aligned}
\|\uv\|_{L^2(0,T;L^2(\Omega))}^2&=
-\int_0^T\int_\Omega\int_0^t\uv(s)\, ds\, \uv_t(t)\, dx\, dt+\int_\Omega \int_0^t\uv(s)\, ds\, \uv(T)\, dx\\
&\leq \Bigl(T \|\uv_t\|_{L^2(0,T;L^2(\Omega))} +\sqrt{T}\|\uv(T)\|_{L^2(\Omega)}\Bigr) \|\uv\|_{L^2(0,T;L^2(\Omega))}\\
\|\uv(T)\|_{L^2(\Omega)}&=\|\int_0^T\uv_t(t)\, dt\|_{L^2(\Omega)}\leq \sqrt{T}\|\uv_t\|_{L^2(0,T;L^2(\Omega))}\,.
\end{aligned}
\]
This, together with the definition of the Neumann Laplacian \eqref{Laplace}, and the use of solutions 
$\zv$, $\vv$ to the auxiliary problems
\begin{equation}\label{auxprob}
\left\{\begin{array}{rcl}
\zv_t-\Delta \zv&=&\vv\mbox{ in }(0,T)\times\Omega\\
\partial_\nu\zv&=&0\mbox{ on }(0,T)\times\partial\Omega\\
\zv(0)&=&0\mbox{ in }\Omega
\end{array}\right.\,, \quad
\left\{\begin{array}{rcl}
-\vv_t-\Delta \vv&=&\mathbf{f}\mbox{ in }(0,T)\times\Omega\\
\partial_\nu\vv&=&0\mbox{ on }(0,T)\times\partial\Omega\\
\vv(T)&=&\mathbf{g}\mbox{ in }\Omega
\end{array}\right.\,,
\end{equation}
allows to derive the identity
\begin{equation}\label{id_adjU}
\begin{aligned}
(\uv,\zv)_\cU 
=& \int_0^T\int_\Omega \Bigl( \nabla\uv\colon\nabla(-\Delta_N\zv)-\uv\cdot\zv_{tt}\Bigr) dx\, dt 
+\int_\Omega \uv(T)\cdot\Bigl(\zv_t(T)-\Delta_N\zv(T)\Bigr)\, dx\\
=& \int_0^T\int_\Omega \Bigl( \nabla\uv\colon\nabla(\vv-\zv_t)-\uv\cdot(\vv_t+\Delta_N\zv_t)\Bigr) dx\, dt 
+\int_\Omega \uv(T)\cdot\vv(T)\, dx\\
=& \int_0^T\int_\Omega \uv\cdot\Bigl(-\Delta_N\vv-\vv_t\Bigr) dx\, dt 
+\int_\Omega \uv(T)\cdot\vv(T)\, dx\\
=&\int_0^T\int_\Omega \uv\cdot \mathbf{f} \, dx\, dt + \int_\Omega \uv(T)\cdot \mathbf{g} \, dx\,,
\end{aligned}
\end{equation}
which will be needed later on for deriving the adjoint.

On $\cW=H^1(0,T;H^1(\Omega;\R^3))^*$ we use the inner product  
\[
\begin{aligned}
(\wv_1,\wv_2)_\cW
&:=\int_0^T\int_\Omega \Bigl( I_1[\nabla(-\Delta_N+\mbox{id})^{-1}\wv_1](t)\colon I_1[\nabla(-\Delta_N+\mbox{id})^{-1}\wv_2](t)\\
&\qquad\qquad+I_1[(-\Delta_N+\mbox{id})^{-1}\wv_1](t)\cdot I_1[(-\Delta_N+\mbox{id})^{-1}\wv_2](t)
\, dx\, dt\,,
\end{aligned}
\]
with the isomorphism $-\Delta_N+\mbox{id}:H^1(\Omega)\to (H^1(\Omega))^*$ and the time integral operators
\[
\begin{aligned}
I_1[w](t):=\int_0^t w(s)\, ds-\frac{1}{T}\int_0^T(T-s)w(s)\, ds\,,\\
I_2[w](t):=-\int_0^t (t-s)w(s)\, ds + \frac{t}{T}\int_0^T(T-s)w(s)\, ds\,,
\end{aligned}
\]
so that $I_2[w]_t(t)=-I_1[w](t)$, $I_1[w]_t(t)=-I_2[w]_{tt}(t)=w(t)$ and $I_2[w](0)=I_2[w](T)=0$, hence 
\[
\int_0^T I_1[w_1](t)\,I_1[w_2](t)\, dt=\int_0^T I_2[w_1](t)\, w_2(t)\, dt,
\]
so that in case $\wv_2\in L^2(0,T;L^2(\Omega;\R^3))$, 
\begin{equation}\label{id_adjW}
\begin{aligned}
(\wv_1,\wv_2)_\cW
&=\int_0^T\int_\Omega \Bigl( I_2[\nabla(-\Delta_N+\mbox{id})^{-1}\wv_1](t)\colon [\nabla(-\Delta_N+\mbox{id})^{-1}\wv_2](t)\\
&\qquad\qquad+I_2[(-\Delta_N+\mbox{id})^{-1}\wv_1](t)\cdot [(-\Delta_N+\mbox{id})^{-1}\wv_2](t)
\, dx\, dt\\
&=\int_0^T\int_\Omega  I_2[(-\Delta_N+\mbox{id})^{-1}\wv_1](t)\cdot\wv_2(t) \, dx\, dt\,.
\end{aligned}
\end{equation}
In case $p>2$ in the assumption on $\hv$, we can set $\cW=H^1(0,T;L^2(\Omega;\R^3))^*$ and use the simpler inner product 
\[
(\wv_1,\wv_2)_\cW
:=\int_0^T\int_\Omega I_1[\wv_1](t)\cdot I_1[\wv_2](t)\, dx\, dt\,,
\]
which in case $\wv_2\in L^2(0,T;L^2(\Omega;\R^3))$ satisfies
\[
(\wv_1,\wv_2)_\cW
=\int_0^T\int_\Omega I_2[\wv_1](t)\cdot \wv_2(t)\, dx\, dt\,.
\]

\subsubsection{Well-definedness of the forward operator}
Indeed it can be verified that $\bF$ maps between the function spaces introduced above, cf. \eqref{U}, \eqref{W}. For the linear (with respect to $\mh$) parts 
$\altil_1 \mh_t$,
$-\Delta_N \mh$, and 
$\int_0^\Tt \int_\Omega \ker_{k\ell}(t,\tau,x)\cdot \mv_t(x,\tau)\, dx\, d\tau$ of $\bF$,
this is obvious and for the nonlinear terms 
$\altil_2 (\mv_0+\mh)\times \mh_t$,
$|\nabla (\mv_0+\mh)|^2(\mv_0+\mh)$, 
$((\mv_0+\mh)\cdot\hv)(\mv_0+\mh)$
we use the following estimates, holding for any $\uv,\wv,\zv\in \cU$.
For the term $\altil_2 (\mv_0+\mh)\times \mh_t$, we estimate
\begin{equation}\label{est_nl1}
\begin{aligned}
\|\uv\times\wv_t\|_{H^1(0,T;H^1(\Omega;\R^3))^*}
&\leq \|\uv\times\wv_t\|_{L^2(0,T;(H^1(\Omega;\R^3))^*)}\\
&\leq C_{H^1\to L^3}^\Omega \|\uv\times\wv_t\|_{L^2(0,T;L^{3/2}(\Omega;\R^3))}\\
&\leq C_{H^1\to L^3}^\Omega \|\uv\|_{C(0,T;L^6(\Omega;\R^3))} \|\wv_t\|_{L^2(0,T;L^2(\Omega;\R^3))}\\
&\leq C_{H^1\to L^3}^\Omega C_{H^1\to L^6}^\Omega\|\uv\|_{C(0,T;H^1(\Omega;\R^3))} \|\wv_t\|_{L^2(0,T;L^2(\Omega;\R^3))}\,,
\end{aligned}
\end{equation}
where we have used duality and continuity of the embeddings $H^1(0,T;H^1(\Omega;\R^3))\hookrightarrow L^2(0,T;H^1(\Omega;\R^3))\hookrightarrow L^2(0,T;L^3(\Omega))$ in the first and second estimate, and H\"older's inequality with exponent $4$ in the third estimate;\\
For the term $|\nabla (\mv_0+\mh)|^2(\mv_0+\mh)$, we use
\begin{equation}\label{est_nl2}
\begin{aligned}
&\|(\nabla\uv\colon\nabla\wv)\zv\|_{H^1(0,T;H^1(\Omega;\R^3))^*}\\
&\leq C_{H^1\to L^\infty}^{(0,T)} \|(\nabla\uv\colon\nabla\wv)\zv\|_{L^1(0,T;(H^1(\Omega;\R^3))^*)}\\
&\leq C_{H^1\to L^\infty}^{(0,T)} C_{H^1\to L^6}^\Omega \|(\nabla\uv\colon\nabla\wv)\zv\|_{L^1(0,T;L^{6/5}(\Omega;\R^3))}\\
&\leq C_{H^1\to L^\infty}^{(0,T)} C_{H^1\to L^6}^\Omega \|\nabla\uv\|_{L^2(0,T;L^6(\Omega;\R^3))}\|\nabla\wv\|_{L^2(0,T;L^6(\Omega;\R^3))} \|\zv\|_{C(0,T;L^2(\Omega;\R^3))}\\
&\leq C_{H^1\to L^\infty}^{(0,T)} (C_{H^1\to L^6}^\Omega;\R^3) \|\uv\|_{L^2(0,T;H^2(\Omega;\R^3))} \|\wv\|_{L^2(0,T;H^2(\Omega;\R^3))} \|\zv\|_{C(0,T;H^1(\Omega;\R^3))}\,,
\end{aligned}
\end{equation}
again using duality and the embeddings $H^1(0,T;H^1(\Omega;\R^3))\hookrightarrow L^\infty(0,T;H^1(\Omega))\hookrightarrow L^\infty(0,T;L^6(\Omega))$;\\
For the term $((\mv_0+\mh)\cdot\hv)(\mv_0+\mh)$, we estimate
\begin{equation}\label{est_nl3}
\begin{aligned}
&\|(\uv\cdot\hv)\zv\|_{H^1(0,T;H^1(\Omega;\R^3))^*}\\
&\leq C_{H^1\to L^6}^\Omega \|(\uv\cdot\hv)\zv\|_{L^2(0,T;L^{6/5}(\Omega;\R^3))}\\
&\leq C_{H^1\to L^6}^\Omega \|\uv\|_{C(0,T;L^6(\Omega;\R^3))} \|\zv\|_{C(0,T;L^6(\Omega;\R^3))} \|\hv\|_{L^2(0,T;L^2(\Omega;\R^3))}\\
&\leq (C_{H^1\to L^6}^\Omega;\R^3) \|\uv\|_{C(0,T;H^1(\Omega;\R^3))} \|\zv\|_{C(0,T;H^1(\Omega;\R^3))} \|\hv\|_{L^2(0,T;L^2(\Omega;\R^3))}
\end{aligned}
\end{equation}
by duality and the embedding $H^1(0,T;H^1(\Omega;\R^3))\hookrightarrow L^2(0,T;L^6(\Omega))$, as well as H\"older's inequality.

In case $p>2$, $\mathbb{F}$ maps into the somewhat stronger space $\cW=H^1(0,T;L^2(\Omega;\R^3))^*$, due to the estimates 
\begin{equation}\label{est_nl1_L2}
\begin{aligned}
\|\uv\times\wv_t\|_{H^1(0,T;L^2(\Omega;\R^3))^*}
&\leq C_{H^1\to L^\infty}^{(0,T)} \|\uv\times\wv_t\|_{L^1(0,T;L^2(\Omega;\R^3))}\\
&\leq C_{H^1\to L^\infty}^{(0,T)} \|\uv\|_{L^2(0,T;L^\infty(\Omega;\R^3))} \|\wv_t\|_{L^2(0,T;L^2(\Omega;\R^3))}\\
&\leq C_{H^1\to L^\infty}^{(0,T)} C_{H^2\to L^\infty}^\Omega\|\uv\|_{L^2(0,T;H^2(\Omega;\R^3))} \|\wv_t\|_{L^2(0,T;L^2(\Omega;\R^3))}\,,
\end{aligned}
\end{equation}
\begin{equation}\label{est_nl2_L2}
\begin{aligned}
&\|(\nabla\uv\colon\nabla\wv)\zv\|_{H^1(0,T;L^2(\Omega;\R^3))^*}\\
&\leq C_{H^1\to L^\infty}^{(0,T)} \|(\nabla\uv\colon\nabla\wv)\zv\|_{L^1(0,T;L^2(\Omega;\R^3))}\\
&\leq C_{H^1\to L^\infty}^{(0,T)} \|\nabla\uv\|_{L^2(0,T;L^6(\Omega;\R^3))}\|\nabla\wv\|_{L^2(0,T;L^6(\Omega;\R^3))} \|\zv\|_{C(0,T;L^6(\Omega;\R^3))}\\
&\leq C_{H^1\to L^\infty}^{(0,T)} (C_{H^1\to L^6}^\Omega;\R^3) \|\uv\|_{L^2(0,T;H^2(\Omega;\R^3))} \|\wv\|_{L^2(0,T;H^2(\Omega;\R^3))} \|\zv\|_{C(0,T;H^1(\Omega;\R^3))}\,,
\end{aligned}
\end{equation}
\begin{equation}\label{est_nl3_L2}
\begin{aligned}
&\|(\uv\cdot\hv)\zv\|_{H^1(0,T;L^2(\Omega;\R^3))^*}\\
&\leq C_{H^1\to L^\infty}^{(0,T)}  \|(\uv\cdot\hv)\zv\|_{L^1(0,T;L^2(\Omega;\R^3))}\\
&\leq C_{H^1\to L^\infty}^{(0,T)}  \|\uv\|_{L^4(0,T;L^{p^{**}}(\Omega;\R^3))} \|\zv\|_{L^4(0,T;L^{p^{**}}(\Omega;\R^3))} \|\hv\|_{L^2(0,T;L^p(\Omega;\R^3))}\\
&\leq C_{H^1\to L^\infty}^{(0,T)} (C_{H^{1/4},L^4}^{(0,T)})^2 (C_{H^{3/2},L^{p^{**}}}^\Omega)^2 \\
& \qquad\qquad \|\uv\|_{H^{1/4}(0,T;H^{3/2}(\Omega;\R^3))} \|\zv\|_{H^{1/4}(0,T;H^{3/2}(\Omega;\R^3))} \|\hv\|_{L^2(0,T;L^p(\Omega;\R^3))}\,,
\end{aligned}
\end{equation}
for $p^{**}=\frac{2p}{p-2}<\infty$, which can be bounded by the $\cU$ norm of $\uv$ and $\zv$, using interpolation with $s=\frac14$ in \eqref{U}.

\subsubsection{Differentiability of the forward operator} 
Formally, the derivative of $\bF$ is given by 
\[
\begin{aligned}
&\bF'(\mh,\altil_1,\altil_2)(\uv,\beta_1,\beta_2)\\
&=\left(\begin{array}{l}
\beta_1 \mh_t - \beta_2 (\mv_0+\mh)\times \mh_t \\
\hspace*{1.5cm}
+\altil_1 \uv_t-\Delta_N \uv - \altil_2 \uv\times \mh_t-\altil_2 (\mv_0+\mh)\times \uv_t
 \\
\hspace*{1.5cm}-2(\nabla (\mv_0+\mh)\colon\nabla\uv)(\mv_0+\mh)- |\nabla (\mv_0+\mh)|^2 \uv
 \\
\hspace*{1.5cm}
+((\mv_0+\mh)\cdot\hv) \uv 
+ (\uv\cdot\hv) (\mv_0+\mh)
\\[2ex]
\Bigl(\int_0^\Tt \int_\Omega \ker_{k\ell}(t,\tau,x) \cdot\uv_t(x,\tau)\, dx\, d\tau\Bigr)_{k=1,\ldots,K\,, \ \ell=1,\ldots,L}
\end{array}\right)\\
&=\left(\begin{array}{ccc}
\frac{\partial\bF_0}{\partial\mh}(\mh,\altil)&\frac{\partial\bF_0}{\partial\altil_1}(\mh,\altil)&\frac{\partial\bF_0}{\partial\altil_2}(\mh,\altil)\\
(\frac{\partial\bF_{k\ell}}{\partial\mh}(\mh,\altil))_{k=1,\ldots,K,\ell=1,\ldots,L}&0&0
\end{array}\right)
\left(\begin{array}{c} \uv \\ \beta_1\\ \beta_2\end{array}\right)
\end{aligned}
\]
where $\frac{\partial\bF_0}{\partial\mh}(\mh,\altil):\cU\to\cW$, $\frac{\partial\bF_0}{\partial\altil_1}(\mh,\altil):\R\to\cW$, $\frac{\partial\bF_0}{\partial\altil_2}(\mh,\altil):\R\to\cW$, $(\frac{\partial\bF_{k\ell}}{\partial\mh}(\mh,\altil))_{k=1,\ldots,K,\ell=1,\ldots,L}:\cU\to L^2(0,T)^{KL}$. 
Fr\'{e}chet differentiability follows from the fact that in 
\[
\bF(\mh+\uv,\altil_1+\beta_1,\altil_2+\beta_2)-\bF(\mh,\altil_1,\altil_2)-\bF'(\mh,\altil_1,\altil_2)(\uv,\beta_1,\beta_2)
\]
all linear terms cancel out and the nonlinear ones are given by (abbreviating $\mv=\mv_0+\mh$)
\[
\begin{aligned}
&(\altil_1+\beta_1)(\mv_t+\uv_t)-\altil_1\mv_t-\altil_1\uv_t-\beta_1 \mv_t 
= \beta_1 \uv_t\\
&(\altil_2+\beta_2)(\mv+\uv)\times(\mv_t+\uv_t) - \altil_2\mv\times\mv_t - \beta_2\mv\times\mv_t - \altil_2\uv\times\mv_t - \altil_2\mv\times\uv_t\\
&\qquad= \altil_2\uv\times\uv_t+\beta_2\mv\times\uv_t+\beta_2\uv\times\mv_t+\beta_2\uv\times\uv_t\\
&|\nabla\mv+\nabla\uv|^2(\mv+\uv) - |\nabla\mv|^2\mv -2 (\nabla\mv\colon\nabla\uv)\mv-|\nabla\mv|^2\uv\\
&\qquad= |\nabla\uv|^2(\mv+\uv)+2(\nabla\mv\colon\nabla\uv)\uv\\
&((\mv+\uv)\cdot\hv)(\mv+\uv)-(\mv\cdot\hv)\mv-(\uv\cdot\hv)\mv-(\mv\cdot\hv)\uv
= (\uv\cdot\hv)\uv\,,
\end{aligned}
\]
hence, using again \eqref{est_nl1}--\eqref{est_nl3}, they can be estimated by some constant multiplied by $\|\uv\|_{\cU}^2+\beta_1^2+\beta_2^2$.

\subsubsection{Adjoints} 

We start with the adjoint of $\frac{\partial\bF_0}{\partial\mh}(\mh,\altil)$. For any $\uv\in\cU$, $\yv\in L^2(0,T;L^2(\Omega))$, we have, using the definition of $-\Delta_N$, i.e., \eqref{Laplace},
\[
\begin{aligned}
&\int_0^T\int_\Omega (\frac{\partial\bF_0}{\partial\mh}(\mh,\altil)\uv)\cdot \yv \,dx\,dt\\
&=
\int_0^T\int_\Omega 
\Bigl(\altil_1 \uv_t\cdot\yv+\nabla \uv\colon\nabla \yv - \altil_2 (\uv\times \mh_t)\cdot\yv-\altil_2 ((\mv_0+\mh)\times \uv_t)\cdot\yv
 \\
&\quad-2(\nabla (\mv_0+\mh)\colon\nabla\uv)\, ((\mv_0+\mh)\cdot\yv)- |\nabla (\mv_0+\mh)|^2 \, (\uv\cdot\yv)
 \\
&\quad+((\mv_0+\mh)\cdot\hv) \,(\uv\cdot\yv) 
+ (\uv\cdot\hv) \,((\mv_0+\mh)\cdot\yv)\Bigr)\,dx\,dt\\
&=
\int_0^T\int_\Omega \uv\cdot
\Bigl(-\altil_1 \yv_t+ (-\Delta\yv) - \altil_2 \mh_t\times\yv+\altil_2 \yv_t\times(\mv_0+\mh)+\altil_2 \yv\times\mh_t
\\&\qquad\qquad\qquad-2((\mv_0+\mh)\cdot\yv)\, (-\Delta_N(\mv_0+\mh))
+2((\nabla (\mv_0+\mh)^T(\nabla \yv))\, (\mv_0+\mh)\\
&\qquad\qquad\qquad +2((\nabla (\mv_0+\mh)^T(\nabla (\mv_0+\mh)))\, \yv
-|\nabla (\mv_0+\mh)|^2 \yv
 \\
&\qquad\qquad\qquad+((\mv_0+\mh)\cdot\hv) \,\yv 
+ ((\mv_0+\mh)\cdot\yv) \, \hv
\Bigr)\,dx\,dt\\
&\quad + \int_\Omega \uv(T)\cdot\Bigl(\altil_1\yv(T) -\altil_2 \yv(T)\times(\mv_0+\mh(T))\Bigr)\, dx\\
&=:\int_0^T\int_\Omega \uv\cdot \mathbf{f}^{\yv}\, dx\, dt + \int_\Omega \uv(T)\cdot \mathbf{g}^{\yv}_T\, dx \,,
\end{aligned}
\]
where we have integrated by parts with respect to time and used the vector identities 
\[
\vec{a}\cdot(\vec{b}\times\vec{c})=\vec{b}\cdot(\vec{c}\times\vec{a})=\vec{c}\cdot(\vec{a}\times\vec{b})\,.
\]
Matching the integrals over $\Omega\times(0,T)$ and $\Omega\times\{T\}$, respectively, and taking into account the homogeneous Neumann boundary conditions implied by the definition of $-\Delta_N$, \eqref{Laplace}, as well as the identities \eqref{id_adjU}, \eqref{id_adjW},
we find that $\frac{\partial\bF_0}{\partial\mh}(\mh,\altil)^*\yv=:\zv$ is the solution of \eqref{auxprob} with $\mathbf{f}=\mathbf{f}^{\yv}$, $\mathbf{g}=\mathbf{g}^{\yv}_T$, 
where in case $\cW=H^1(0,T;H^1(\Omega;\R^3))^*$, $\yv=I_2[\widetilde{y}]$, with  $\widetilde{y}(t)$ solving
\[
\left\{\begin{array}{rcl}
-\Delta \widetilde{y}(t)+\widetilde{y}(t)&=&\wv(t)\mbox{ in }\Omega\\
\partial_\nu\widetilde{y}&=&0\mbox{ on }\partial\Omega
\end{array}\right.
\]
for each $t\in(0,T)$, 
or in case $\cW=H^1(0,T;L^2(\Omega;\R^3))^*$, just $\yv=I_2[\wv]$.

With the same $\yv$, after pointwise projection onto the mutually orthogonal vectors $\mh_t(x,t)$ and $(\mv_0(x)+\mh(x,t))\times\mh_t(x,t)$ and integration over space and time, we also get the adjoints of $\frac{\partial\bF_0}{\partial\altil_1}(\mh,\altil)$, $\frac{\partial\bF_0}{\partial\altil_2}(\mh,\altil)$
\begin{align*}
\frac{\partial\bF_0}{\partial\altil_1}(\mh,\altil)^*\wv &= \int_0^T\int_\Omega \mh_t\cdot \yv\, dx\, dt\,, \\
\frac{\partial\bF_0}{\partial\altil_2}(\mh,\altil)^*\wv &= -\int_0^T\int_\Omega ((\mv_0+\mh)\times\mh_t)\cdot \yv\, dx\, dt\,. 
\end{align*}
Finally, the fact that for $\uv\in\cU$, $y\in L^2(0,T)^{KL}$
\begin{equation}\label{adjA21}
\begin{aligned}
&\left(\Big(\frac{\partial\bF_{k\ell}}{\partial\mh}(\mh,\altil)\Big)_{k=1,\ldots,K,\ell=1,\ldots,L}\uv,y\right)_{L^2(0,T)^{KL}}\\
&=\sum_{k=1}^K\sum_{\ell=1}^L \int_0^T \left(\Big(\frac{\partial\bF_{k\ell}}{\partial\mh}(\mh,\altil)\Big)_{k=1,\ldots,K,\ell=1,\ldots,L}\uv\right)_{k\ell}(t) y_{k\ell}(t)\, dt\\
&=\sum_{k=1}^K\sum_{\ell=1}^L \int_0^T \int_0^\Tt \int_\Omega \ker_{k\ell}(t,\tau,x)\cdot \uv_t(x,\tau)\, dx\, d\tau y_{k\ell}(t)\, dt\\
&=\sum_{k=1}^K\sum_{\ell=1}^L \int_0^T \Bigl(
-\int_0^\Tt \int_\Omega \frac{\partial}{\partial\tau}\ker_{k\ell}(t,\tau,x)\cdot \uv(x,\tau)\, dx\, d\tau \\
& \hspace{3cm} + \int_\Omega \ker_{k\ell}(t,\Tt,x)\cdot \uv(x,\Tt)\, dx
\Bigr) y_{k\ell}(t)\, dt\,,
\end{aligned}
\end{equation}
where we have integrated by parts with respect to time, implies that due to \eqref{id_adjU},
\begin{displaymath}
 (\frac{\partial\bF_{k\ell}}{\partial\mh}(\mh,\altil))_{k=1,\ldots,K,\ell=1,\ldots,L}^*y=\zv
\end{displaymath}
is obtained by solving another auxiliary problem \eqref{auxprob} with 
\begin{equation}\label{fgA21}
\begin{split}
\mathbf{f}(x,\tau)&=-\int_0^T\sum_{k=1}^K\sum_{\ell=1}^L \frac{\partial}{\partial\tau}\ker_{k\ell}(t,\tau,x) y_{k\ell}(t)\ dt,\\
\mathbf{g}(x)&=\int_0^T\sum_{k=1}^K\sum_{\ell=1}^L \ker_{k\ell}(t,\Tt,x) y_{k\ell}(t)\ dt\,.
\end{split}
\end{equation}
\begin{remark}
In case of a Landweber-Kaczmarz method iterating cyclically over the equations defined by $\bF_0,\bF_{k\ell}$, $k=1,...,K$, $\ell=1,...,L$, adjoints of derivatives of $\bF_0$ remain unchanged while adjoints of $\frac{\partial\bF_{k\ell}}{\partial\mh}(\mh,\altil))_{k=1,\ldots,K,\ell=1,\ldots,L}$ are defined as in \eqref{adjA21}, \eqref{fgA21} by just skipping the sums over $k$ and $\ell$ there.
\end{remark}

%% file: diffadj_red.tex
\subsection{Reduced formulation}\label{sec:red}
We now consider the formulation \eqref{ip_red} with $F$ defined by \eqref{red-F}, \eqref{defS}, \eqref{defK}.
Due to the estimate 
\begin{align*}
\|\mathcal{K}_{k\ell}\m_t\|^2_{L^2(0,T)}\leq T\|\widetilde{a}_{\ell}\|^2_{L^2(0,T)} \|c_k\pr\|^2_{L^2(\Om2R)}\|\m\|^2_{H^1(0,T;L^2(\Om2R))}\,,
\end{align*}
if $\widetilde{a}_{\ell}\in L^2(0,T), c_k\pr\in L^2(\Om2R)$ we can choose the state space in the reduced setting as
\begin{align}\label{red-spaceU}
\cUtil=H^1(0,T;L^2(\Om2R)),
\end{align}
which is different from the one in the all-at-once setting.

\subsubsection{Adjoint equation}
From \eqref{red-F} the derivative of the forward operation takes the form
\begin{align}\label{red-diffF}
F'(\altil)\beta= \cK\u_t,
\end{align}
where $\u$ solves the linearized LLG equation
\begin{alignat*}{3}
&\altil_1\u_t - \altil_2\m\times\u_t - \altil_2\u\times\m_t - \Delta\u -2(\nabla\u:\nabla\m)\m\\
&\qquad+\u(-|\nabla\m|^2+(\m\cdot \h)) + (\u\cdot \h)\m &&\\
&=-\beta_1\m_t+\beta_2\m\times\m_t \qquad &&\text{in } (0,T)\times\Omega\\
& \partial_\nu\u=0 && \text{on } (0,T)\times\partial\Omega\\
& \u(0)=0 && \text{in } \Omega,
\end{alignat*}
and $\m$ is the solution to \eqref{llg3}-\eqref{llg3_abs}.
This equation can be obtained by formally taking directional derivatives (in the direction of $\u$) in all terms of the LLG equation \eqref{llg3}--\eqref{llg3_abs}, or alternatively by subtracting the defining boundary value problems for $S(\m+\epsilon\u)$ and $S(\m)$, dividing by $\epsilon$ and then letting $\epsilon$ tend to zero.

The Hilbert space adjoint
\begin{align*}
F'(\altil)^*:L^2(0,T)^{KL}\rightarrow\R^2
\end{align*}
of $F'(\altil)$ satisfies, for each $z\in L^2(0,T)^{KL}$,
{\allowdisplaybreaks
\begin{align}\label{red-adjoint0}
&(F'(\altil)^*z,\beta)_{\R^2} \nonumber\\
&=(z,F'(\altil)\beta)_{L^2(0,T)^{KL}} \nonumber\\
&=\sum_{k=1}^K\sum_{\ell=1}^L\int_0^T z_{k\ell}(t) \int_0^T\int_\Omega (-\mu_0)\widetilde{a}_\ell(t-\tau) c_k(x)\pr(x)\cdot\u_\tau(\tau,x) dx\,d\tau\,dt \nonumber\\
&=\sum_{k=1}^K\sum_{\ell=1}^L\int_0^T z_{k\ell}(t)\bigg(-\int_0^T\int_\Omega (-\mu_0)\cdot (-1)\widetilde{a}_{\ell\,t}(t-\tau) c_k(x)\pr(x)\cdot\u(\tau,x)\,dx\,d\tau \nonumber\\
&\hspace{2.5cm} +\int_\Omega (-\mu_0)\widetilde{a}_{\ell}(t-T)c_k(x)\pr(x)\cdot\u(T,x)\,dx \bigg)dt  \nonumber\\
&= \int_0^T\int_\Omega\u(\tau,x)\cdot \sum_{k=1}^K\sum_{\ell=1}^L \bigg( \int_0^T (-\mu_0)\widetilde{a}_{\ell\,t}(t-\tau)z_{k\ell}(t)\,dt \bigg)\, c_k(x)\pr(x)\,dx\,d\tau \nonumber\\
&\hspace{2.5cm} + \int_\Omega\u(T,x)\cdot \sum_{k=1}^K\sum_{\ell=1}^L \bigg(\int_0^T (-\mu_0)\widetilde{a}_{\ell}(t)z_{k\ell}(t)\,dt \bigg)\, c_k(x)\pr(x)\,dx  \nonumber\\
&=: (\u,\Ktil z)_{L^2(0,T;L^2(\Om2R))} + (\u(T),\KtilT z)_{L^2(\Om2R)}
\end{align}
}
as the transfer function $\tilde{a}$ is periodic with period $T$, and the continuous embedding $H(0,T)\hookrightarrow C[0,T]$ allows us to evaluate $\u(t=T)$.\\

Observing
{\allowdisplaybreaks
\begin{alignat*}{3}
&\int_0^T\int_\Omega -\altil_1\p_t\cdot\u\,dx\,dt
=\int_0^T\int_\Omega\altil_1\u_t\cdot\p \,dx - \int_\Omega\altil_1\p(T)\cdot\u(T)\,dx \,, \\
&\int_0^T\int_\Omega -\altil_2(\m\times\p)_t\cdot\u \,dx\,dt\\
&=\int_0^T\int_\Omega -\altil_2 (\m\times\u_t)\cdot\p\,dx\,dt -\int_\Omega\altil_2 (\m\times\p)(T)\cdot\u(T)\,dx \,,\\
&\int_0^T\int_\Omega\altil_2(\p\times\m_t)\cdot\u \,dx\,dt
=\int_0^T\int_\Omega -\altil_2(\u\times\m_t)\cdot\p \,dx\,dt \,,\\
&\int_0^T\int_\Omega -\Delta\p\cdot\u \,dx\,dt
=\int_0^T\int_\Omega -\p\cdot\Delta\u \,dx\,dt - \int_0^T\int_{\partial\Omega} \partial_\nu\p\cdot\u \,dx\,dt \,,\\
&\int_0^T\int_\Omega \p(-|\nabla\m|^2+(\m\cdot\h))\cdot\u \,dx\,dt
=\int_0^T\int_\Omega \left(\u(-|\nabla\m|^2+(\m\cdot \h))\right)\cdot\p \,dx\,dt \,,\\
&\int_0^T\int_\Omega \left(\p\cdot\m\right)\h\cdot\u\,dx\,dt
=\int_0^T\int_\Omega (\u\cdot\h)\m\cdot\p\,dx\,dt \,,\\
&\int_0^T\int_\Omega 2(\m\cdot\p)\Delta\m\cdot\u \,dx\,dt\\
&=-\int_0^T\int_\Omega 2(\nabla\m:\nabla\u)(\m\cdot\p) \,dx\,dt\\
&\qquad\qquad\qquad +2\int_0^T\int_\Omega -\u\cdot((\nabla\m)^\top\nabla\m)\p -\u\cdot((\nabla\m)^\top\nabla\p)\m \,dx\,dt \,,
\end{alignat*}
}
we see that, if $\p$ solves the adjoint equation
\begin{alignat}{3}
&-\altil_1\p_t - \altil_2\m\times\p_t - 2\altil_2\m_t\times\p - \Delta\p \nonumber\\
&\quad +2\left((\nabla\m)^\top\nabla\m\right)\p + 2\left((\nabla\m)^\top\nabla\p\right)\m \nonumber\\
&\quad+(-|\nabla\m|^2+(\m\cdot \h))\p + (\m\cdot\p)(\h + 2\Delta\m)=\Ktil z  \qquad&&\text{in } (0,T)\times\Omega \label{red-adjoint-eq0-1}\\
& \partial_\nu\p=0 && \text{on } (0,T)\times\partial\Omega \label{red-adjoint-eq0-2}\\
& \altil_1\p(T)+\altil_2(\m\times\p)(T)=\KtilT z && \text{in } \Omega \label{red-adjoint-eq0-3}
\end{alignat}
then with \eqref{red-adjoint0}, we have
\begin{align*}
(F'(\altil)^*z,\beta)_{\R^2}&=(\u,\Ktil z)_{L^2(0,T;L^2(\Om2R))} + (\u(T),\KtilT z)_{L(\Om2R)}\\
&=\int_0^T\int_\Omega(-\beta_1\m_t+\beta_2\m\times\m_t)\cdot\p \,dx\,dt \nonumber\\
&=(\beta_1,\beta_2)\cdot\left( \int_0^T\int_\Omega -\m_t\cdot\p \,dx\,dt , \int_0^T\int_\Omega(\m\times\m_t)\cdot\p \,dx\,dt\right),
\end{align*}
which implies the Hilbert space adjoint  $F'(\altil)^*:\cY\rightarrow\R^2$
\begin{align} \label{red-adjoint}
F'(\altil)^*z=\left( \int_0^T\int_\Omega -\m_t\cdot\p \,dx\,dt , \int_0^T\int_\Omega(\m\times\m_t)\cdot\p \,dx\,dt\right),
\end{align}
provided that the adjoint state $\p$ exists and belongs to a sufficiently smooth space (see Subsection \ref{sec-red-adjoint-solvability} below).

The final condition \eqref{red-adjoint-eq0-3} is equivalent to
\begin{align*}
\begin{pmatrix}
\altil_1 & -\altil_2\m_3(T) & \altil_2\m_2(T)\\
\altil_2\m_3(T) & \altil_1 & -\altil_2\m_1(T)\\
-\altil_2\m_2(T) & \altil_2\m_1(T) & \altil_1
\end{pmatrix}\p(T)=:M^{\altil}_T\p(T)=\KtilT z,
\end{align*}
where $\m_i(T), i=1,2,3,$ denotes the i-th component of $\m(T)$. The matrix $M^{\altil}_T$ with $\det(M^{\altil}_T)=|\altil_1(\altil_1^2+\altil_2^2)|$ is invertible if $\altil_1 > 0$, which matches the condition for existence of the solution to the LLG equation. Hence, we are able to rewrite the adjoint equation in the form
\begin{alignat}{3}
&-\altil_1\p_t - \altil_2\m\times\p_t - 2\altil_2\m_t\times\p - \Delta\p \nonumber\\
&\quad +2\left((\nabla\m)^\top\nabla\m\right)\p + 2\left((\nabla\m)^\top\nabla\p\right)\m \nonumber\\
&\quad+(-|\nabla\m|^2+(\m\cdot \h))\p + (\m\cdot\p)(\h + 2\Delta\m)=\Ktil z  \qquad&&\text{in } (0,T)\times\Omega \label{red-adjoint-eq-1}\\
& \partial_\nu\p=0 && \text{on } (0,T)\times\partial\Omega \label{red-adjoint-eq-2}\\
& \p(T)=(M_T^{\altil})^{-1}\KtilT z && \text{in } \Omega. \label{red-adjoint-eq-3}
\end{alignat}

\begin{remark}
Formula \eqref{red-adjoint} inspires a Kaczmarz scheme relying on restricting the observation operator to time subintervals for every fixed $k, \ell$, namely, we segment $(0,T)$ into several subintervals $(t^j, t^{j+1})$ with the break points $0=t^0<\ldots<t^{n-1}=T$ and
\begin{align}
F^j_{k\ell}:\mathcal{D}(F)(\subseteq\cX)\rightarrow\cY^j, \qquad \altil \mapsto y^j
:=\mathcal{K}_{k\ell} \frac{\partial}{\partial t} S(\altil)|_{(t^j, t^{j+1})} 
\end{align}
with
\begin{align}
\cY^j=L^2(t^j, t^{j+1})^{KL} \qquad j=0\ldots n-1,
\end{align}
hence
\begin{align}
y^j_{k\ell}(t)=\int_{t^j}^{t^{j+1}}\int_\Omega -\mu_0\widetilde{a}_{\ell}(t-\tau) c_k(x)\pr(x)\cdot\m_\tau(x,\tau) dx d\tau.
\end{align}
Here we distinguish between the superscript $j$ for the time subinterval index and subscripts $k, \ell$ for the index of different receive coils and concentrations. 

For $z^j\in\cY^j$, 
\begin{align*}
&(\Ktil z^j)(x,t)=\sum_{k=1}^K\sum_{\ell=1}^L-\mu_0c_k(x)\pr(x)\int_{t^j}^{t^{j+1}} \widetilde{a}_{\ell\,\tau}(\tau-t)z^j_{k\ell}(\tau)\,d\tau \qquad t\in(0,T) \,,\\
&(\KtilT z^j)(x)=\sum_{k=1}^K\sum_{\ell=1}^L-\mu_0c_k(x)\pr(x)\int_{t^j}^{t^{j+1}} \widetilde{a}_{\ell}(\tau)z^j_{k\ell}(\tau)\,d\tau
\end{align*}
yield the same Hilbert space adjoint $F^{j'}(\altil)^*:\cY^j\rightarrow\R^2$ as in \eqref{red-adjoint}, and the adjoint state $\pp^{z^j}$ still needs to be solved on the whole time line $[0,T]$ with
\begin{alignat}{3}
&-\altil_1\pp^{z^j}_t - \altil_2\m\times\pp^{z^j}_t - 2\altil_2\m_t\times\pp^{z^j} - \Delta\pp^{z^j} \nonumber\\
&\quad +2\left((\nabla\m)^\top\nabla\m\right)\pp^{z^j} + 2\left((\nabla\m)^\top\nabla\pp^{z^j}\right)\m \nonumber\\
&\quad+(-|\nabla\m|^2+(\m\cdot \h))\pp^{z^j} + (\m\cdot\pp^{z^j})(\h + 2\Delta\m)=\Ktil z^j  \quad&&\text{in } (0,T)\times\Omega \label{red-adjoint-eqLK-1}\\
& \partial_\nu\pp^{z^j}=0 && \text{on } (0,T)\times\partial\Omega \label{red-adjoint-eqLK-2}\\
& \pp^{z^j}(T)=(M_T^{\altil})^{-1}\KtilT z^j && \text{in } \Omega. \label{red-adjoint-eqLK-3}
\end{alignat}

Besides this, the conventional Kaczmarz method resulting from the collection of observation operators $\cK_{k\ell}$ with $k=1\ldots K, \ell=1\ldots L$ as in \eqref{forwardoperator_kl} is always applicable, where
\begin{align}
F_{k\ell}:\mathcal{D}(F)(\subseteq\cX)\rightarrow\cY_{k\ell}, \qquad \altil \mapsto y_{k\ell}:=\mathcal{K}_{k\ell}\frac{\partial}{\partial t}(S(\altil))
\end{align}
with
\begin{align}
\cY_{k\ell}=L^2(0,T)\qquad k=1\ldots K, \ \ell=1\ldots L
\end{align}
Thus ${F'_{k\ell}(\altil)}^*$ can be seen as \eqref{red-adjoint}, where the adjoint state $\p_{k\ell}$ solves \eqref{red-adjoint-eq-1}-\eqref{red-adjoint-eq-3} with corresponding data
\begin{align*}
&\Ktil_{k\ell}z(x,t)=-\mu_0c_k(x)\pr(x)\int_0^T \widetilde{a}_{\ell\,\tau}(\tau-t)z(\tau)\,d\tau \qquad t\in(0,T) \,,\\
&\Ktil_{T\,k\ell}z(x)=-\mu_0c_k(x)\pr(x)\int_0^T \widetilde{a}_{\ell}(\tau)z(\tau)\,d\tau
\end{align*}
for each $z\in\cY_{k\ell}$.
\end{remark}

\subsubsection{Solvability of the adjoint equation} \label{sec-red-adjoint-solvability}
First of all, we derive a bound for $\p$. To begin with, we set $\tau=T-t$ to convert \eqref{red-adjoint-eq-1}-\eqref{red-adjoint-eq-3} into an initial boundary value problem. Then we test \eqref{red-adjoint-eq-1} with $\p_t$ and obtain the identities and estimates
{\allowdisplaybreaks
\begin{alignat*}{3}
&\int_\Omega \altil_1\p_t(t)\cdot\p_t(t) \,dx=\altil_1\|\p_t(t)\|^2_{L^2(\Om2R)} \,,\\
&\int_\Omega \altil_2(\m(t)\times\p_t(t))\cdot\p_t(t) \,dx=0 \,,\\
&\int_\Omega \altil_2(\m_t(t)\times\p(t))\cdot\p_t(t) \,dx
\leq |\altil_2| \|\m_t(t)\|_{L^3(\Om2R)}\|\p(t)\|_{L^6(\Om2R)}\|\p_t(t)\|_{L^2(\Om2R)} \,,\\
&\int_\Omega -\Delta\p(t)\cdot\p_t(t) \,dx=\frac{1}{2}\frac{d}{dt}\|\nabla\p(t)\|^2_{L^2(\Om2R)} \,,\\
&\int_\Omega \left(((\nabla\m(t))^\top\nabla\m(t))\p(t)\right)\cdot\p_t(t) \,dx\\
&\, \leq (C^{\Omega}_{H^1\rightarrow L^6})^2\|\nabla\m\|^2_{L^\infty(0,T;H^1(\Om2R))}\|\p(t)\|_{L^6(\Om2R)}\|\p_t(t)\|_{L^2(\Om2R)} \,,\\
&\int_\Omega \left(((\nabla\m(t))^\top\nabla\p(t))\m(t)\right)\cdot\p_t(t) \,dx \\
&\,\leq C^\Omega_{H^2\rightarrow L^\infty}\|\nabla\m(t)\|_{H^2(\Om2R)}\|\nabla\p(t)\|_{L^2(\Om2R)}\|\p_t(t)\|_{L^2(\Om2R)} \,,\\
&\int_\Omega (-|\nabla\m(t)|^2+(\m(t)\cdot \h))\p(t)\cdot\p_t(t) \,dx\\
&\, \leq \left((C^{\Omega}_{H^1\rightarrow L^6})^2 \|\nabla\m\|^2_{L^\infty(0,T;H^1(\Om2R))}+\|\h(t)\|_{L^3(\Om2R)} \right)\|\p(t)\|_{L^6(\Om2R)}\|\p_t(t)\|_{L^2(\Om2R)} \,,\\
&\int_\Omega \left(\m(t)\cdot\p(t)\right)\h(t)\cdot\p_t(t) \,dx
\leq \|\h(t)\|_{L^3(\Om2R)}\|\p(t)\|_{L^6(\Om2R)}\|\p_t(t)\|_{L^2(\Om2R)} \,,\\
&\int_\Omega (\m(t)\cdot\p(t))\Delta\m(t)\cdot\p_t(t) \,dx\\
&\,\leq C^\Omega_{H^1\rightarrow L^3}\|\Delta\m(t)\|_{H^1(\Om2R))}\|\p(t)\|_{L^6(\Om2R)}\|\p_t(t)\|_{L^2(\Om2R)} \,,\\
&\int_\Omega \Ktil z(t)\cdot\p_t(t) \,dx
\leq \|\Ktil z(t)\|_{L^2(\Om2R)}\|\p_t(t)\|_{L^2(\Om2R)}\,.
\end{alignat*}
}
Above, we employ the fact that the solution $\m$ to the LLG equation has $|\m|=1$ and the continuity of the embeddings $H^1(\Om2R)\hookrightarrow L^6(\Om2R)\hookrightarrow L^3(\Om2R), H^2(\Om2R)\hookrightarrow L^\infty(\Om2R)$ through the constants $C^\Omega_{H^1\rightarrow L^6}, C^\Omega_{H^1\rightarrow L^3}$ and $C^\Omega_{H^2\rightarrow L^\infty}$, respectively.\\
Employing Young's inequality we deduce, for each $t\leq T$ and $\epsilon>0$ sufficiently small,
\begin{align} \label{red-p-bound-0}
&\frac{1}{2}\frac{d}{dt}\|\nabla\p(t)\|^2_{L^2(\Om2R)}+(\altil_1-\epsilon)\|\p_t(t)\|^2_{L^2(\Om2R)} \nonumber\\
&\leq \bigg[\left(\|\nabla\m\|^4_{L^\infty(0,T;H^1(\Om2R))}+\|\nabla\m(t)\|^2_{H^2(\Om2R)}+ \|\m_t(t)\|^2_{L^3(\Om2R)}+\|\h(t)\|^2_{L^3(\Om2R)}\right) \nonumber\\
&\qquad\qquad\quad .\|\p(t)\|^2_{H^1(\Om2R)} +\|\Ktil z(t)\|^2_{L^2(\Om2R)}\bigg]\frac{C}{4\epsilon}.
\end{align}
The generic constant $C$ might take different values whenever it appears.

To have the full $H^1-$norm on the left hand side of this estimate, we apply the transformation $\tilde{\p}(t)=e^t\p(t)$, which yields $\tilde{\p}_t(t)=e^{t}(\p(t)+\p_t(t))$. After testing by $\p_t$, the term $\int_\Omega \p(t)\cdot\p_t(t) \,dx=\frac{1}{2}\frac{d}{dt}\|\p(t)\|^2_{L^2(\Om2R)}$ will contribute to $\frac{1}{2}\frac{d}{dt}\|\nabla\p(t)\|^2_{L^2(\Om2R)}$ forming the full $H^1-$norm on the left hand side. Alternatively, one can add $\p$ to both sides of \eqref{red-adjoint-eq-1} and evaluate the right hand side with $\int_\Omega \p(t)\cdot\p_t(t) \,dx\leq \frac{1}{4\epsilon}\|\p(t)\|^2_{H^1(\Om2R)}+\epsilon\|\p_t(t)\|^2_{L^2(\Om2R)}$.

Integrating over $(0,t)$, we get
\begin{align*} 
&\frac{1}{2}\|\p(t)\|^2_{H^1(\Om2R)}+(\altil_1-\epsilon)\|\p_t\|^2_{L^2(0,t;L^2(\Om2R))} \nonumber\\
&\leq \frac{C}{4\epsilon} \bigg[\int_0^t \Big(\|\nabla\m\|^4_{L^\infty(0,T;H^1(\Om2R))}+\|\nabla\m(\tau)\|^2_{H^2(\Om2R)}+ \|\m_t(\tau)\|^2_{L^3(\Om2R)}\\
&\qquad\qquad\qquad +\|\h(\tau)\|^2_{L^3(\Om2R)}\Big).\|\p(\tau)\|^2_{H^1(\Om2R)}\,d\tau \\
&\qquad\qquad\qquad +\|\Ktil z\|^2_{L^2(0,T;L^2(\Om2R))} + \|(M_T^{\altil})^{-1}\KtilT z\|^2_{H^1(\Om2R)}\bigg]
\end{align*}
with the evaluation for the terms $\|\Ktil z\|_{L^2(0,T;L^2(\Om2R))}$ and $\|(M_T^{\altil})^{-1}\KtilT z\|^2_{H^1(\Om2R)}$ (not causing any misunderstanding, we omit here the subscripts $k,\ell$ for indices of concentrations and coil sensitivities)
\begin{alignat*}{3}
&\|\Ktil z(t)\|_{L^2(\Om2R)}^2\leq C\|c\textbf{p}^R\|^2_{L^2(\Om2R)}\|\tilde{a}\|^2_{H^1(0,T)}\|z\|^2_{L^2(0,T)}\leq C^{\tilde{a},c,\textbf{p}^R}\|z\|^2_{L^2(0,T)} \,,\\[1ex]
&|(M_T^{\altil})^{-1}\KtilT z\|^2_{H^1(\Om2R)}\\
&\leq C^{\altil}\|z\|^2_{L^2(0,T)}\|\tilde{a}\|^2_{L^2(0,T)}\\
&\qquad\quad.\left( \|c\textbf{p}^R\|^2_{H^1(\Om2R)}+\|c\textbf{p}\m_i(T)\|^2_{H^1(\Om2R)}+\|c\textbf{p}^R\m_j(T)\m_k(T)\|^2_{H^1(\Om2R)} \right)\\
&\leq C^{\altil_0,\rho, \tilde{a}}\|z\|^2_{L^2(0,T)}\left(\|c\textbf{p}^R\|^2_{H^1(\Om2R)}+\|c\textbf{p}^R\|^2_{L^6(\Om2R)}\|\nabla\m(T)\|^2_{L^3(\Om2R)} \right)\\
&\leq  C^{\tilde{a}}\|z\|^2_{L^2(0,T)} \\
&\qquad\quad.\left( \|c\textbf{p}^R\|^2_{H^1(\Om2R)}+(C^{\Omega}_{H^1\rightarrow L^6}C^{\Omega}_{H^1\rightarrow L^3})^2\|c\textbf{p}^R\|^2_{H^1(\Om2R)}\|\nabla\m\|^2_{L^\infty(0,T;H^1(\Om2R))}  \right)\\
&\leq C^{\tilde{a},c,\textbf{p}^R}\|z\|^2_{L^2(0,T)}\|\nabla\m\|^2_{L^\infty(0,T;H^1(\Om2R))}
\end{alignat*}
with some $i,j,k=1,2,3$. This estimate holds for $c\textbf{p}^R\in H^1(\Om2R)$ and thus requires some smoothness of the concentration $c$, while the coil sensitivity $\textbf{p}^R$ is usually smooth in practice.

Then applying Gr\"onwall's inequality yields 
\begin{align*}
&\|\p\|_{L^\infty(0,T;H^1(\Om2R))}\\
&\leq C\exp\Big(\|\nabla\m\|^2_{L^\infty(0,T;H^1(\Om2R))}+\|\nabla\m\|_{L^2(0,T;H^2(\Om2R))}+\|\m_t\|_{L^2(0,T;L^3(\Om2R))}\\
&\qquad\quad + \|\h\|_{L^2(0,T;L^3(\Om2R))}\Big)
 .\big(\|\Ktil z\|_{L^2(0,T;L^2(\Om2R))}+\|(M_T^{\altil})^{-1}\KtilT z\|_{H^1(\Om2R)} \big)\\
&\leq C^{\tilde{a},c,\textbf{p}^R}\Big(\|\nabla\m\|_{L^\infty(0,T;H^1(\Om2R))\cap L^2(0,T;H^2(\Om2R))}, \|\m_t\|_{L^2(0,T;L^3(\Om2R))}\\
&\qquad\quad ,\|\h\|_{L^2(0,T;L^3(\Om2R))}\Big).\|z\|_{L^2(0,T)}.
\end{align*}
Integrating \eqref{red-p-bound-0} on $(0,T)$, we also get
\begin{align*}
&\|\p_t\|_{L^2(0,T;L^2(\Om2R))}\\
&\leq C^{\tilde{a},c,\textbf{p}^R}\Big(\|\nabla\m\|_{L^\infty(0,T;H^1(\Om2R))\cap L^2(0,T;H^2(\Om2R))}, \|\m_t\|_{L^2(0,T;L^3(\Om2R))})\\
&\qquad\qquad\quad ,\|\h\|_{L^2(0,T;L^3(\Om2R))}\Big) .\|z\|_{L^2(0,T)}.
\end{align*}
Altogether, we obtain
\begin{align} \label{red-p-bound}
&\|\p\|_{L^\infty(0,T;H^1(\Om2R))}+\|\p_t\|_{L^2(0,T;L^2(\Om2R))}\nonumber \\
&\leq C^{\tilde{a},c,\textbf{p}^R}\Big(\|\nabla\m\|_{L^\infty(0,T;H^1(\Om2R))\cap L^2(0,T;H^2(\Om2R))}, \|\m_t\|_{L^2(0,T;L^3(\Om2R))} \nonumber\\
&\qquad\qquad\quad ,\|\h\|_{L^2(0,T;L^3(\Om2R))}\Big) .\|z\|_{L^2(0,T)}. 
\end{align}
This result applied to the Galerkin approximation implies existence of the solution to the adjoint equation. Uniqueness also follows from  \eqref{red-p-bound}.

\subsubsection{Regularity of the solution to the LLG equation}
In \eqref{red-p-bound}, first of all we need the solution $\m\in L^\infty(0,T;H^2(\Om2R))$\\$\cap L^2(0,T;H^3(\Om2R))$ to the LLG equation. This can be obtained from the regularity result in \cite[Lemma 2.3]{bgmh93} for $\m_0\in{H^2(\Om2R)}$ with small $\|\nabla\m_0\|_{L^2(\Om2R)}$. The remaining task is verifying that the estimate still holds in case the external field $\h$ is present, i.e., the right hand side of \eqref{llg3} contains the additional term $\Pr\h$.

Following the lines of the proof in \cite[Lemma 2.3]{bgmh93}, we take the second spatial derivative of $\Pr\h$, then test it by $\Delta\m$ such that
\begin{alignat*}{4}
&\int_\Omega  \Delta\h(t)\cdot\Delta\m(t) \,dx\\
&\leq \begin{cases}
\|\Delta\h(t)\|_{L^2(\Om2R)}\|\Delta\m(t)\|_{L^2(\Om2R)}  &\text{if } \h\in L^2(0,T;H^2(\Om2R))\\
\|\nabla\h(t)\|_{L^2(\Om2R)}\|\nabla^3\m(t)\|_{L^2(\Om2R)}  &\text{if } \h\in L^2(0,T;H^1(\Om2R)),\, \partial_\nu\h=0 \text{ on } \partial\Omega\\
\end{cases} \,,\\[1ex]
&\int_\Omega  \Delta((\m(t)\cdot\h(t))\m(t))\cdot\Delta\m(t) \,dx\\
&\leq \begin{cases}
C\|\h(t)\|_{H^2(\Om2R)}\left(1+6\|\nabla\m(t)\|_{H^1(\Om2R)}+2\|\nabla\m(t)\|_{H^2(\Om2R)}\|\nabla\m\|_{L^\infty(0,T;L^2(\Om2R))}\right)\\
\hspace{1cm} .\|\Delta\m(t)\|_{L^2(\Om2R)} 
\hspace{1.5cm}\, \text{ if } \h\in L^2(0,T;H^2(\Om2R))\\
C\|\h(t)\|_{H^1(\Om2R)}\left( 1+2\|\nabla\m(t)\|_{L^3(\Om2R)} \right)\|\nabla^3\m(t)\|_{L^2(\Om2R)}\\ \hspace{5cm}\text{if } \h\in L^2(0,T;H^1(\Om2R)),\, \partial_\nu\h=0 \text{ on } \partial\Omega\\
\end{cases}
\end{alignat*}
with $C$ just depending on the constants in the embeddings $H^1(\Om2R)\hookrightarrow L^6(\Om2R)\hookrightarrow L^3(\Om2R)$. Then we can proceed similarly to the proof of \cite[Lemma 2.3]{bgmh93} by applying Young's inequality, Gronwall's inequality and time integration to arrive at
\begin{align}\label{red-GHlem-h}
&\|\nabla \m\|_{L^\infty(0,T;H^1(\Om2R))\cap L^2(0,T;H^2(\Om2R))} \nonumber\\
&\hspace{3cm}\leq \left(\|\nabla \m_0\|_{H^1(\Om2R)}+\|\h\|\right)C(\|\nabla \m_0\|_{H^1(\Om2R)},\|\h\|),
\end{align}
where $\|\h\|$ is evaluated in $L^2(0,T;H^1(\Om2R))$ or $L^2(0,T;H^2(\Om2R))$ as in the two cases mentioned above.\\

It remains to prove $\m_t\in L^2(0,T;H^1(\Om2R))\hookrightarrow L^2(0,T;L^3(\Om2R))$ to validate \eqref{red-p-bound}. For this purpose, instead of working with \eqref{llg3} we test \eqref{llg2} by $-\Delta\m_t$ and obtain
{\allowdisplaybreaks
\begin{align*}
&\int_\Omega \m_t(t)\cdot(-\Delta\m_t(t)) \,dx=\|\nabla\m_t(t)\|^2_{L^2(\Om2R)} \,,\\
&\int_\Omega -\alpha_1 \Delta\m(t)\cdot(-\Delta\m_t(t)) \,dx=\frac{\alpha_1}{2}\frac{d}{dt}\|\Delta\m(t)\|^2_{L^2(\Om2R)} \,,\\
&\int_\Omega -\alpha_1 |\nabla\m(t)|^2\m(t)\cdot(-\Delta\m_t(t)) \,dx
= -\alpha_1\int_\Omega \nabla\left(|\nabla\m(t)|^2\m(t)\right):\nabla\m_t(t) \,dx\\
&\,\leq \alpha_1\Big(2C^\Omega_{H^1\rightarrow L^6}C^\Omega_{H^1\rightarrow L^3}\|\nabla\m\|_{L^\infty(0,T;H^1(\Om2R))}\|\Delta\m(t)\|_{H^1(\Om2R)}\\
&\qquad\qquad\quad + (C^{\Omega}_{H^1\rightarrow L^6})^3\|\nabla\m\|^3_{L^\infty(0,T;H^1(\Om2R))}  \Big).\|\nabla\m_t(t)\|_{L^2(\Om2R)} \,,\\
&\int_\Omega -\alpha_1 (\h(t)-(\m(t)\cdot\h(t))\m(t))\cdot(-\Delta\m_t(t)) \,dx\\
& =-\alpha_1\int_\Omega \nabla(\h(t)-(\m(t)\cdot\h(t))\m(t)):\nabla\m_t(t) \,dx\\
&\,\leq 2\alpha_1\Big(\|\nabla\h(t)\|_{L^2(\Om2R)} \\
&\qquad\qquad\quad + C^\Omega_{H^1\rightarrow L^6}\|\h(t)\|_{L^3(\Om2R)}\|\nabla\m\|_{L^\infty(0,T;H^1(\Om2R))} \Big). \|\nabla\m_t(t)\|_{L^2(\Om2R)} \,,\\
&\int_\Omega -\alpha_2(\m(t)\times\Delta\m(t))\cdot(-\Delta\m_t(t)) \,dx =\int_\Omega -\alpha_2\nabla(\m(t)\times\Delta\m(t)):\nabla\m_t(t) \,dx\\
&\,\leq |\alpha_2|\Big(C^\Omega_{H^1\rightarrow L^6}C^\Omega_{H^1\rightarrow L^3}\|\nabla\m\|_{L^\infty(0,T;H^1(\Om2R))}\|\Delta\m(t)\|_{H^1(\Om2R)}\\
&\qquad\qquad\quad + \|\nabla^3\m(t)\|_{L^2(\Om2R)} \Big).\|\nabla\m_t(t)\|_{L^2(\Om2R)} \,,\\
&\int_\Omega -\alpha_2(\m(t)\times\h(t))\cdot(-\Delta\m_t(t)) \,dx=\int_\Omega -\alpha_2\nabla(\m(t)\times\h(t)):(\nabla\m_t(t)) \,dx\\
&\,\leq |\alpha_2|\Big(C^\Omega_{H^1\rightarrow L^6}\|\h(t)\|_{L^3(\Om2R)}\|\nabla\m\|_{L^\infty(0,T;H^1(\Om2R))}\\
&\qquad\qquad\quad + \|\nabla\h(t)\|_{L^2(\Om2R)} \Big).\|\nabla\m_t(t)\|_{L^2(\Om2R)} \,.
\end{align*}
}
Integrating over $(0,T)$ then employing H\"older's inequality,  Young's inequality and \eqref{red-GHlem-h}, it follows that
\begin{align}
&(1-\epsilon)\|\nabla\m_t\|_{L^2(0,T;L^2(\Om2R))}\nonumber\\
&\leq \frac{C}{4\epsilon} \Big(\|\nabla\m\|_{L^\infty(0,T;H^1(\Om2R))}\|\nabla\m\|_{L^2(0,T;H^2(\Om2R))}+\|\nabla\m\|^3_{L^\infty(0,T;H^1(\Om2R))} \nonumber\\
&\qquad\qquad +\|\nabla\m\|_{L^2(0,T;H^2(\Om2R))}  + \|\h\|_{L^2(0,T;H^1(\Om2R))}\|\nabla\m\|_{L^\infty(0,T;H^1(\Om2R))}  \nonumber\\
&\qquad\qquad +\|\h\|_{L^2(0,T;H^1(\Om2R))}\Big) \nonumber\\
&\leq \left(\|\nabla \m_0\|_{H^1(\Om2R)}+\|\h\|\right)C(\|\nabla \m_0\|_{H^1(\Om2R)},\|\h\|).
\end{align}
Also $\|\m_t\|_{L^2(0,T;L^2(\Om2R))}<C\left(\|\nabla\m_0\|_{L^2(\Om2R)}+\|\h\|_{L^2(0,T;L^2(\Om2R))}\right)$ according to \cite{mkap06} with taking into account the presence of $\h$, we arrive at
\begin{align}\label{Mregularity}
\|\m_t\|_{L^2(0,T;H^1(\Om2R))} \leq \left(\|\nabla \m_0\|_{H^1(\Om2R)}+\|\h\|\right)C(\|\nabla \m_0\|_{H^1(\Om2R)},\|\h\|),
\end{align}
where $\|\h\|$ is evaluated in $L^2(0,T;H^1(\Om2R))$ or $L^2(0,T;H^2(\Om2R))$.

In conclusion, the fact that $\m\in L^\infty(0,T;H^2(\Om2R))\cap L^2(0,T;H^3(\Om2R))\cap H^1(0,T;H^1(\Om2R))$ for $\m_0\in H^2(\Om2R)$ with small $\|\nabla\m_0\|_{L^2(\Om2R)}$, and\\
 $\h\in L^2(0,T;H^1(\Om2R)), \partial_\nu\h=0$ on $\partial\Omega$ or $\h\in L^2(0,T;H^2(\Om2R))$ guarantee unique existence of the adjoint state $\p\in L^\infty(0,T;H^1(\Om2R))\cap H^1(0,T;L^2(\Om2R))$. And this regularity of $\p$ ensures the adjoint $F'(\altil)^*$ in \eqref{red-adjoint} to be well-defined.

\begin{remark}\text{ }
\begin{enumerate}[label=$\bullet$]
\item
The LLG equation \eqref{llg3}-\eqref{llg3_abs} is uniquely solvable for $\altil_1>0$ and arbitrary $\altil_2$. Therefore, the regularization problem should be locally solved within the ball $\cB_\rho(\altil^0)$ of center $\altil^0$ with $\altil^0_1>0$ and radius $\rho<\altil_1^0$.
\item
\cite[Lemma 2.3]{bgmh93} requires smallness $\|\nabla\m_0\|_{L^2(\Om2R)}\leq\lambda$, and this smallness depends on $\altil$ through the relation $C^I\left(\lambda^2+2\lambda+\frac{\altil_2}{\altil_1}\lambda \right)<1$ with $C^I$ depending on the constants in the interpolation inequalities.
\end{enumerate}
\medskip
Altogether, we arrive at
\begin{align}
\mathcal{D}(F)=\left\{\altil=(\altil_1,\altil_2)\in\cB_\rho(\altil^0): 0<\altil^0_1,\rho<\altil_1^0, C^I\left(\lambda^2+2\lambda+\frac{\altil_2}{\altil_1}\lambda \right)<1 \right\}.
\end{align}
\end{remark}

\subsubsection{Differentiability of the forward operator}
\label{sec-red-diff}
Since the observation operator $\mathcal{K}$ is linear, differentiability of $F$ is just the question of differentiability of $S$.

Let us rewrite the LLG equation \eqref{llg3} in the following form
\begin{align*}
\tilde{g}(\altil,\m)-\Delta\m = \tilde{f}(\m)
\end{align*}
and denote
\begin{align*}
\vtil^\epsilon:=\frac{S(\altil+\epsilon\beta)-S(\altil)}{\epsilon}-\u
=:\frac{\n-\m}{\epsilon}-\u=:\ve-\u.
\end{align*}
Considering the system of equations
\begin{alignat*}{4}
&\tilde{g}(\altil+\epsilon\beta,\n) &&-\Delta\n &&=\tilde{f}(\n),\\
&\tilde{g}(\altil,\m) &&-\Delta\m &&=\tilde{f}(\m),\\
&\tilde{g}'_\m(\altil,\m)\u + \tilde{g}'_{\altil}(\altil,\m)\beta &&-\Delta\u && =\tilde{f}'_\m(\m)\u,
\end{alignat*}
with the same boundary and initial data for each, we see that $\vtil^\epsilon$ solves
\begin{alignat}{3}
&\tilde{g}'_\m(\altil,\m)\vtil^\epsilon -\Delta\vtil^\epsilon -\tilde{f}'_\m(\m)\vtil^\epsilon  \nonumber\\
&= \frac{\tilde{f}(\n)-\tilde{f}(\m)}{\epsilon}-\tilde{f}'_\m(\m)\ve - \frac{\tilde{g}(\altil+\epsilon\beta,\n)-\tilde{g}(\altil,\m)}{\epsilon}&& \label{red-vtil-1}\\
&\qquad  + \tilde{g}'_\m(\altil,\m)\ve + \tilde{g}'_{\altil}(\altil,\m)\beta &&\text{in } (0,T)\times\Omega \nonumber\\
& \partial_\nu\vtil^\epsilon=0 && \text {on } [0,T]\times\partial\Omega \label{red-vtil-2}\\
& \vtil^\epsilon(0)=0 && \text{in } \Omega,\label{red-vtil-3}
\end{alignat}
explicitly
\begin{alignat}{3}
&\altil_1\vtil^\epsilon_t - \altil_2\m\times\vtil^\epsilon_t - \altil_2\vtil^\epsilon\times\m_t - \Delta\vtil^\epsilon \nonumber\\
&\quad -2(\nabla\vtil^\epsilon:\nabla\m)\m +\vtil^\epsilon(-|\nabla\m|^2+(\m\cdot \h)) + (\vtil^\epsilon\cdot \h)\m \nonumber\\
&=\frac{1}{\epsilon}\left( |\nabla\n|^2\n+\Prn\h -|\nabla\m|^2\m-\Pr\h \right) \label{red-vtil-4}\\
&\quad-2(\nabla\ve:\nabla\m)\m+\ve(-|\nabla\m|^2+(\m\cdot \h)) + (\ve\cdot \h)\m \nonumber\\
&\quad -\frac{1}{\epsilon}\Big( (\altil_1+\epsilon\beta_1)\n_t - (\altil_2+\epsilon\beta_2)\n\times\n_t  -\altil_1\m_t + \altil_2\m \times &&\m_t\Big) \nonumber\\
&\quad +\altil_1\vet - \altil_2\m\times\vet - \altil_2\ve\times\m_t \nonumber\\
&\quad +\beta_1\m_t-\beta_2\m\times\m_t \qquad &&\text{in } (0,T)\times\Omega  \nonumber\\
& \partial_\nu\vtil^\epsilon=0 &&\text {on } [0,T]\times\partial\Omega \label{red-vtil-5}\\
& \vtil^\epsilon(0)=0 &&\text{in } \Omega.\label{red-vtil-6}
\end{alignat}
Observing the similarity of \eqref{red-vtil-4}-\eqref{red-vtil-6} to the adjoint equation \eqref{red-adjoint-eq-1}-\eqref{red-adjoint-eq-3} with $\vtil^\epsilon$ in place of $\p$ and denoting by $\be$ the right-hand side of \eqref{red-vtil-1} or \eqref{red-vtil-4}, one can evaluate $\|\vtil^\epsilon\|$ using the same technique as in Section \ref{sec-red-adjoint-solvability}. By this way, one achieves, for each $\epsilon\in[0,\bar{\epsilon}]$,
\begin{align*}
\|\vtil^\epsilon\|_{L^\infty(0,T;H^1(\Om2R))\cap H^1(0,T;L^2(\Om2R))}\leq C\|\be\|_{L^2(0,T;L^2(\Om2R))}
\end{align*}
with $\be\in L^2(0,T;L^2(\Om2R))$ also by analogously estimating and employing $\m,\n\in L^\infty(0,T;H^2(\Om2R))\cap L^2(0,T;H^3(\Om2R))\cap H^1(0,T;H^1(\Om2R))$. We note that the constant $C$ here is independent of $\epsilon$.\\

Next letting $\cV:={L^\infty(0,T;H^1(\Om2R))\cap H^1(0,T;L^2(\Om2R))}$, we have
{\allowdisplaybreaks
\begin{alignat*}{3}
&\|\be\|_{L^2(0,T;L^2(\Om2R))}=\bigg\lVert \frac{\tilde{f}(\n)-\tilde{f}(\m)}{\epsilon}-\tilde{f}'_\m(\m)\ve - \frac{\tilde{g}(\altil+\epsilon\beta,\n)-\tilde{g}(\altil,\m)}{\epsilon}\\
&\hspace{4cm}+ \tilde{g}'_\m(\altil,\m)\ve + \tilde{g}'_{\altil}(\altil,\m)\beta \bigg\lVert_{L^2(0,T;L^2(\Om2R))}\\
&\leq \bigg\lVert \int_0^1  \Big( \big(\tilde{f}'_\m(\m+\lambda\epsilon\ve)-\tilde{f}'_\m(\m) \big)\ve
- \big( \tilde{g}'_\m(\altil+\lambda\epsilon\beta,\m+\lambda\epsilon\ve) - \tilde{g}'_\m(\altil,\m) \big)\ve\\
&\hspace{2.5cm} - \big( \tilde{g}'_{\altil}(\altil+\lambda\epsilon\beta,\m+\lambda\epsilon\ve) - \tilde{g}'_{\altil}(\altil,\m) \big)\beta \Big)\,d\lambda  \bigg\lVert_{L^2(0,T;L^2(\Om2R))}\\
&\leq 2\sup_{\substack{\lambda\in[0,1]\\ \epsilon\in[0,\bar{\epsilon}]}}\Big( \|\tilde{f}'_\m(\m+\lambda\epsilon\ve)\|_{\cV\rightarrow L^2(0,T;L^2(\Om2R))}\|\ve\|_{\cV}\\
&\hspace{2.5cm} + \|\tilde{g}'_\m(\altil+\lambda\epsilon\beta,\m+\lambda\epsilon\ve)\|_{\cV\rightarrow L^2(0,T;L^2(\Om2R))}\|\ve\|_{\cV} \\
&\hspace{2.5cm} + \|\tilde{g}'_{\altil}(\altil+\lambda\epsilon\beta,\m+\lambda\epsilon\ve)\|_{\R^2\rightarrow L^2(0,T;L^2(\Om2R))}|\beta| \Big).
\end{alignat*}
}
In order to prove uniform boundedness of the derivatives of $\tilde{f}, \tilde{g}$ w.r.t $\lambda, \epsilon$ in the above estimate, we again proceed in a similar manner as in Section \ref{sec-red-adjoint-solvability} since the space for $\p$ in Section \ref{sec-red-adjoint-solvability} (c.f. \eqref{red-GHlem-h}) coincides with $\cV$ here and by the fact that 
\begin{align}
\max\{\|\m\|,\|\n\|\}
&\leq \max\Big\{\frac{1}{\altil_1},\frac{1}{\altil_1+\epsilon\beta_1}\Big\} C\left(\|\m_0\|_{H^2(\Om2R))},\|\h\|_{L^2(0,T;H^2(\Om2R))}\right) \nonumber\\
&\leq\frac{C}{{\altil^0_1}-\rho}
\end{align}
for $\m,\n\in L^\infty(0,T;H^2(\Om2R))\cap L^2(0,T;H^3(\Om2R))\cap H^1(0,T;H^1(\Om2R))$. \\
If $\partial_\nu\h=0$ on $\partial\Omega$, we just need the $\|.\|_{L^2(0,T;H^1(\Om2R))}$-norm for $\h$ as claimed in \eqref{red-GHlem-h}. This estimate holds for any $\epsilon\in[0,\bar{\epsilon}]$, and the constant $C$ is independent of $\epsilon$.\\

To accomplish uniform boundedness for $\|\be\|_{L^2(0,T;L^2(\Om2R))}$, we need to show that $\|\ve\|_{\cV}$ is also uniformly bounded w.r.t $\epsilon$. It is seen from
\begin{alignat*}{4}
&\tilde{g}(\altil+\epsilon\beta,\n) &&-\Delta\n &&=\tilde{f}(\n),\\
&\tilde{g}(\altil,\m) &&-\Delta\m &&=\tilde{f}(\m)
\end{alignat*}
that $\ve$ solves
\begin{alignat}{3}
&\int_0^1 \tilde{g}'_\m(\altil+\lambda\epsilon\beta,\m+\lambda\epsilon\ve)\ve +  \tilde{g}'_{\altil}(\altil+\lambda\epsilon\beta,\m+\lambda\epsilon\ve)\beta \,d\lambda && - \Delta\ve \nonumber\\
&\qquad= \int_0^1 \tilde{f}'_\m(\m+\lambda\epsilon\ve)\ve\, d\lambda 
 \qquad&&\text{in } (0,T)\times\Omega \label{red-ve-1}\\
& \partial_\nu\ve=0 \qquad &&\text {on } [0,T]\times\partial\Omega \label{red-ve-2}\\
& \ve(0)=0 && \text{in } \Omega\label{red-ve-3}.
\end{alignat}
Noting that $\M:=\m+\lambda\epsilon\ve=\lambda\n+(1-\lambda)\m$  has $\|\M\|\leq \frac{C}{{\altil^0_1}-\rho}$ for all $\lambda\in[0,1]$  with $C$ being independent of $\epsilon$, and $\tilde{g}$ is first order in $\altil$, we can rewrite \eqref{red-ve-1} into the linear equation
\begin{align}
\tilde{G}(\altil+\lambda\epsilon\beta,\M)\ve  -\Delta\ve + \tilde{F}(\M)\ve=\tilde{B}(\M)\beta.
\end{align}
Following the lines of the proof in Section \ref{sec-red-adjoint-solvability}, boundedness of the terms $-\Delta, \tilde{F}(\M)$, $\tilde{B}(\M)$ are straightforward, while the main term in $\tilde{G}(\altil+\lambda\epsilon\beta,\M)$ producing the single square norm of $\vet$, after being tested by $\vet$ is
\begin{align*}
\int_0^1 (\altil_1+\lambda\epsilon\beta_1)\int_\Omega \vet(t)\cdot\vet(t)\,dx\,d\lambda &= \|\vet(t)\|^2_{L^2(\Om2R)}\left(\altil_1+\frac{\epsilon\beta_1}{2}\right)\\
&\geq \|\vet(t)\|^2_{L^2(\Om2R)}(\altil^0_1-\rho).
\end{align*}
According to this, one gets, for all $\epsilon\in[0,\bar{\epsilon}]$,
\begin{align} \label{red-ve-bound}
\|\ve\|_\cV \leq C|\beta|\|\tilde{B}(\M)\|_{\R^2\rightarrow L^2(0,T;L(\Om2R))}\leq |\beta|C
\end{align}
with $C$ depending only on $\m_0, \h, \altil^0,\rho$.\\

Since $\be\rightarrow 0$ pointwise and $\|\be\|_{L^2(0,T;L^2(\Om2R))}\leq C$ for all $\epsilon\in[0,\bar{\epsilon}]$, applying Lebesgue's Dominated Convergence Theorem yields convergence of $\|\be\|_{L^2(0,T;L^2(\Om2R))}$, thus of $\|\vtil^\epsilon\|_\cV$, to zero. Fr\'echet differentiability of the forward operator in the reduced setting is therefore proved.

%% file: conclusion.tex
\section{Conclusion} \label{sec_concl}

In this contribution we outlined a mathematical model of MPI taking into account relaxation effects, which led us to the LLG equation describing the behavior of the magnetic material inside the particles on a microscale level. For calibrating the MPI device it is necessary to compute the system function, which mathematically can be interpreted as an inverse parameter identification problem for an initial boundary value problem based on the LLG equation. To this end we deduced a detailed analysis of the forward model, i.e., the operator mapping the coefficients to the solution of the PDE as well as of the underlying inverse problem. The inverse problem itself was investigated in an all-at-once and a reduced approach. The analysis includes representations of the respective adjoint operators and Fr\'{e}chet derivatives. These results are necessary for a subsequent numerical computation of the system function in a robust manner, which will be subject of future research. Even beyond this, the analysis might be useful for the development of solution methods for other inverse problems that are connected to the LLG equation.

%% file: MPI_LLG_arxiv.bbl
\begin{thebibliography}{10}

\bibitem{borcea02}
{\sc L.~Borcea}, {\em {Electrical impedance tomography}}, Inverse Problems, 18
  (2002), pp.~R99--R136.

\bibitem{bsffgppr14}
{\sc F.~Bruckner, D.~Suess, M.~Feischl, T.~F{\"u}hrer, P.~Goldenits, M.~Page,
  D.~Praetorius, and M.~Ruggeri}, {\em {Multiscale modeling in micromagnetics:
  Existence of solutions and numerical integration}}, Mathematical Models and
  Methods in Applied Sciences, 24 (2014), pp.~2627--2662.

\bibitem{ck13}
{\sc D.~Colton and R.~Kress}, {\em Inverse Acoustic and Electromagnetic
  Scattering Theory}, Springer New York, 2013.

\bibitem{cgc12}
{\sc L.~R. Croft, P.~W. Goodwill, and S.~M. Conolly}, {\em Relaxation in
  x-space magnetic particle imaging}, IEEE transactions on medical imaging, 31
  (2012), pp.~2335--2342.

\bibitem{cg11imm}
{\sc B.~D. Cullity and C.~D. Graham}, {\em Introduction to {M}agnetic
  {M}aterials}, John Wiley \& Sons, 2011.

\bibitem{demt2}
{\sc W.~Demtroeder}, {\em Experimentalphysik 2}, Springer Berlin Heidelberg,
  2013.

\bibitem{dkps15}
{\sc T.~Dunst, M.~Klein, A.~Prohl, and A.~Sch{\"a}fer}, {\em Optimal control in
  evolutionary micromagnetism}, IMA Journal of Numerical Analysis, 35 (2015),
  pp.~1342--1380.

\bibitem{tlg04}
{\sc T.~L. Gilbert}, {\em A phenomenological theory of damping in ferromagnetic
  materials}, IEEE Transactions on Magnetics, 40 (2004), pp.~3443--3449.

\bibitem{bgjw05}
{\sc B.~Gleich and J.~Weizenecker}, {\em {Tomographic imaging using the
  nonlinear response of magnetic particles}}, Nature, 435 (2005),
  pp.~1214--1217.

\bibitem{bgmh93}
{\sc B.~Guo and M.-C. Hong}, {\em {The Landau-Lifshitz equation of the
  ferromagnetic spin chain and harmonic maps}}, Calculus of Variations and
  Partial Differential Equations, 1 (1993), pp.~311--334.

\bibitem{HKLS07}
{\sc M.~Haltmeier, R.~Kowar, A.~Leitao, and O.~Scherzer}, {\em Kaczmarz methods
  for regularizing nonlinear ill-posed equations {II}: Applications}, Inverse
  Problems and Imaging, 1 (2007), pp.~507--523.

\bibitem{HLS07}
{\sc M.~Haltmeier, A.~Leitao, and O.~Scherzer}, {\em Kaczmarz methods for
  regularizing nonlinear ill-posed equations {I}: convergence analysis},
  Inverse Problems and Imaging, 1 (2007), pp.~289--298.

\bibitem{HaNeSc95}
{\sc M.~Hanke, A.~Neubauer, and O.~Scherzer}, {\em A convergence analysis of
  the {L}andweber iteration for nonlinear ill-posed problems}, Numer.\ Math.,
  72 (1995), pp.~21--37.

\bibitem{HINSHAW;LENT:83}
{\sc W.S. Hinshaw and A.H. Lent}, {\em An introduction to {N}{M}{R} imaging:
  {F}rom the {B}loch equation to the imaging equation}, Proc. IEEE, 71 (1983),
  pp.~338--350.

\bibitem{bk16}
{\sc B.~Kaltenbacher}, {\em Regularization based on all-at-once formulations
  for inverse problems}, SIAM Journal of Numerical Analysis, 54 (2016),
  pp.~2594--2618.

\bibitem{bk17}
{\sc B.~Kaltenbacher}, {\em All-at-once versus reduced iterative methods for
  time dependent inverse problems}, Inverse Problems, 33 (2017), p.~064002.

\bibitem{kns08}
{\sc B.~Kaltenbacher, A.~Neubauer, and O.~Scherzer}, {\em Iterative
  Regularization Methods for Nonlinear Ill-Posed Problems}, De Gruyter, 2008.

\bibitem{kss09}
{\sc B.~Kaltenbacher, F.~Sch{\"o}pfer, and T.~Schuster}, {\em {Iterative
  methods for nonlinear ill-posed problems in {B}anach spaces: convergence and
  applications to parameter identification problems}}, Inverse Problems, 25
  (2009), p.~065003.

\bibitem{akar16}
{\sc A.~Kirsch and A.~Rieder}, {\em {Inverse problems for abstract evolution
  equations with applications in electrodynamics and elasticity}}, Inverse
  Problems, 32 (2016), p.~085001.

\bibitem{tk18}
{\sc T.~Kluth}, {\em Mathematical models for magnetic particle imaging},
  Inverse Problems, 34 (2018), p.~083001.

\bibitem{tkpm17}
{\sc T.~Kluth and P.~Maass}, {\em Model uncertainty in magnetic particle
  imaging: Nonlinear problem formulation and model-based sparse
  reconstruction}, International Journal on Magnetic Particle Imaging, 3
  (2017).

\bibitem{tktb12}
{\sc T.~Knopp and T.~M. Buzug}, {\em Magnetic Particle Imaging: an Introduction
  to Imaging Principles and Scanner Instrumentation}, Springer Berlin
  Heidelberg, 2012.

\bibitem{tkngmm17}
{\sc T.~Knopp, N.~Gdaniec, and M.~M{\"o}ddel}, {\em Magnetic particle imaging:
  from proof of principle to preclinical applications}, Physics in Medicine \&
  Biology, 62 (2017), p.~R124.

\bibitem{mkap06}
{\sc M.~Kruz{\'i}k and A.~Prohl}, {\em {Recent Developments in the Modeling,
  Analysis, and Numerics of Ferromagnetism}}, SIAM Review, 48 (2006),
  pp.~439--483.

\bibitem{llel92}
{\sc L.~Landau and E.~Lifshitz}, {\em {3 - On the theory of the dispersion of
  magnetic permeability in ferromagnetic bodies Reprinted from Physikalische
  Zeitschrift der Sowjetunion 8, Part 2, 153, 1935.}}, in {Perspectives in
  Theoretical Physics}, L.P. PITAEVSKI, ed., Pergamon, Amsterdam, 1992,
  pp.~51--65.

\bibitem{Landweber}
{\sc L.~Landweber}, {\em {An iteration formula for Fredholm integral equations
  of the first kind}}, American Journal of Mathematics, 73 (1951),
  pp.~615--624.

\bibitem{natterer86}
{\sc F.~Natterer}, {\em The Mathematics of Computerized Tomography},
  Vieweg+Teubner Verlag, 1986.

\bibitem{natterer_buch_2}
{\sc F.~Natterer and F.~W{\"u}bbeling}, {\em Mathematical {M}ethods in {I}mage
  {R}econstruction}, SIAM, Philadelphia, 2001.

\bibitem{ttnn19}
{\sc T.~T.~N. Nguyen}, {\em {Landweber{\textendash}Kaczmarz for parameter
  identification in time-dependent inverse problems: all-at-once versus reduced
  version}}, Inverse Problems, 35 (2019), p.~035009.

\bibitem{drjb14}
{\sc D.~B. Reeves and J.~B. Weaver}, {\em Approaches for modeling magnetic
  nanoparticle dynamics}, Critical Reviews in Biomedical Engineering, 42
  (2014).

\bibitem{rieder99}
{\sc A.~Rieder}, {\em {On the regularization of nonlinear ill-posed problems
  via inexact {N}ewton iterations}}, Inverse Problems, 15 (1999), pp.~309--327.

\bibitem{Roubicek}
{\sc T.~Roub{\'\i}{\v{c}}ek}, {\em Nonlinear Partial Differential Equations
  with Applications}, International Series of Numerical Mathematics, Springer
  Basel, 2013.

\bibitem{SKHK12}
{\sc Th. Schuster, B.~Kaltenbacher, B.~Hofmann, and K.S. Kazimierski}, {\em {
  Regularization Methods in Banach Spaces}}, de Gruyter, Berlin, New York,
  2012.
\newblock Radon Series on Computational and Applied Mathematics.

\bibitem{SHEPP;VARDI:82}
{\sc L.A. Shepp and Y.~Vardi}, {\em Maximum likelihood reconstruction for
  emission tomography}, IEEE Trans. Med. Imag., 1 (1982), pp.~113--122.

\bibitem{ws16}
{\sc A.~Wald and T.~Schuster}, {\em Sequential subspace optimization for
  nonlinear inverse problems}, Journal of Inverse and Ill-posed Problems, 25
  (2016), pp.~99--117.

\bibitem{awts17}
\leavevmode\vrule height 2pt depth -1.6pt width 23pt, {\em Tomographic
  terahertz imaging using sequential subspace optimization}, New Trends in
  Parameter Identification, Springer Series Trends in Mathematics,  (2018).

\end{thebibliography}
